\documentclass[11pt]{article}

\usepackage[english]{babel}
\usepackage{AVDoc_Disc}

\let\OLDthebibliography\thebibliography 
\renewcommand\thebibliography[1]{ 
	\OLDthebibliography{#1} 
	\setlength{\parskip}{0pt} 
	\setlength{\itemsep}{5pt plus 0.3ex} 
}

\begin{document}
\title{
Exponential quasi-ergodicity \\
for processes with discontinuous trajectories
 }
\author{Aur\'{e}lien Velleret \thanks{Universit\'{e} Paris-Saclay, INRAE, MaIAGE, Domaine de Vilvert, 78350 Jouy-en-Josas, France; E-mail:
		\texttt{aurelien.velleret@nsup.org}
	}}

\maketitle
\subsection*{Abstract} 

This paper tackles the issue
of establishing a lower-bound 
on the asymptotic ratio of survival probabilities
between two different initial conditions,
asymptotically in time 
for a given Markov process with extinction.
Such a comparison is a crucial
step in recent techniques 
for proving  exponential convergence 
to a quasi-stationary distribution.
We introduce a weak form of the Harnack inequality
as the essential ingredient
for such a comparison.
This property is actually a consequence of the convergence  property
we intend to prove.
Its complexity appears as the price to pay
for the level of flexibility required by our applications.
We show in our illustrations 
how simply and efficiently 
it can be used nonetheless.
As illustrations,
we consider two continuous-time processes on $\bR^d$ 
that do not satisfy  the classical Harnack inequalities, 
even in a local version.
The first one is a piecewise deterministic process
while the second is a pure jump process with restrictions
on the directions of its jumps.
\\

\textit{\textbf{Keywords: }} 
continuous-time and continuous-space Markov process;\; 
jumps;\;
quasi-stationary distribution;\; 
survival capacity;\; 
Q-process;\; 
Harris recurrence 
\\

\textit{\textbf{Mathematics Subject Classification:} 37A30, 28D10, 60J25 and 92D15.}

\section{Introduction}
\label{D_sec_intro}
\setcounter{eq}{0}

\subsection{General presentation}

	This work is concerned 
with the long time behavior of
quite general strong Markov processes,
conditionally upon the fact that this process
has not been absorbed in some "cemetery state"
(i.e. that it is not "extinct").
The eventual interest is on the analogous of stationary distributions when such a conditioning
is taken into account, namely quasi-stationary distributions (QSD).

In the aftermath of recent works by Champagnat and Villemonais, 
notably in \cite{ChQSD},
we are interested in highlighting 
key properties that ensure convergence results at exponential rate
towards QSD. 
While the approach was initiated in the framework of a convergence in total variation
that is uniform over the initial condition, 
how to deal with heterogeneity in the initial condition has already been the concern of further studies
(\cite{CVly2, AV_QSD, BCGM19}).

These works are inspired by the Harris recurrence techniques 
that is exploited for the proof of convergence towards a stationary distribution 
for conservative semi-groups (cf Section 2 in \cite{KM03} or \cite{MT93}, notably chapter 15, for more details). 
The core of these techniques is a Doeblin minorization condition
 where the density of the marginal law is lower-bounded  
uniformly over some initial conditions. 
A Lyapunov criterion is then exploited to deal with the heterogeneity in the initial conditions.
While refinements of these two properties can be identified in \cite{CVly2, AV_QSD, BCGM19},
another key property is introduced for the generalization to non-conservative semi-groups 
that compare survival between different initial conditions.

This kind of property has a different nature as the others in that 
it is expressed asymptotically in large time.
In this paper, as a representative of those properties,
we focus on the following
property stated in Equation \eqref{D_BdSv}.
It is more simply expressed 
in that it is stated uniformly over the initial condition 
and the denominator is taken as a single reference measure $\alc\in \M_1(\cX)$.
\begin{equation}
	\limSInf{t}\; \sup_{x\in \cX} \; \dfrac{\PR_x(t<\ext)}{\PR_{\alc}(t<\ext)} 
	< \infty,
	\label{D_BdSv}
\end{equation}
In the expression, $\ext$ denotes the extinction time
and the initial condition of the process is written in subscript.
Although $\alc$ is meant to relate to small sets,
no property on $\alc$ should be needed for the proof of  \eqref{D_BdSv}.
\eqref{D_BdSv} is the archetype of the comparison properties of asymptotic survival probability, 
that are introduced in the abstract.
\\

The derivation of a property like \eqref{D_BdSv}
is well understood in the cases of discrete-space processes
(cf the Proof of Theorem 4.1 p14 in \cite{ChQSD}),
and processes that satisfy either
two-sided estimates (cf. Section 2 in \cite{CCV18}),
gradient estimates (cf. Section 3 in \cite{CCV18})
or the Harnack inequality
(step 4 of the proof given in Section 4 of \cite{ChpLyap}, or more generally stated in Subsection 4.2.3 in \cite{AV_QSD}).
In the case of many processes involving jumps, and notably in multi-dimensional settings,
none of these three estimates seem to be applicable.
It seems difficult for the other techniques that we know of (notably Lemma 5 in \cite{DV16}) 
to determine to what extent they could be generalised.

This is why we worked at a weakened form of the Harnack inequality with more flexibility,
(that we call "almost perfect harvest", see Subsection \ref{sec_prop}),
that would still suffice to obtain \eqref{D_BdSv},
when we exploit the other key properties requested for exponential convergence towards QSD.
The aim is to capture a much larger range of processes.
In practice, 
the statement of our main result (see Theorem \ref{D_Th:AF})
relies on the other key properties proposed in \cite{AV_QSD}.

This choice is not simply due to the fact 
that we were more involved in the proofs of \cite{AV_QSD},
but also because the fact that the approach is trajectorial
gives more insight into the estimates
and because the proofs in \cite{AV_QSD} were designed 
with the current article already in mind.
Potential generalizations 
are nevertheless discussed in Subsection \ref{D_sec_gen}.

As a direct consequence, 
we present in Subsection \ref{D_sec_EC}
the new set of key conditions
that is proved sufficient  
for the general results of \cite{AV_QSD} to be deduced.
It implies not only 
the existence and uniqueness 
of the QSD, 
but also several results of exponential convergence.
\\

The paper is organized as follows.
In the next Subsection \ref{D_sec_iApp}, 
we present two illustrative models 
that contribute to motivate 
the new conditions introduced in this article. 
With these two practical issues in mind, 
a concise presentation of the results 
can more readily be addressed in Subsection \ref{sec_spec}. 
There, we first specify the convergence objective, 
then how the contribution of \cite{AV_QSD} 
leads to proving such convergence results 
and finally how the present contribution 
eases the verification of one of the main conditions of \cite[Section 2]{AV_QSD}.
In Subsection~\ref{D_sec_lit}, 
we then place our method 
in the context of existing results on quasi-stationarity.

In complement to this concise overview of our results,
a more detailed description is provided in Section \ref{D_sec_MR}.
The focus is then on the main contribution of the article,
namely the proposed key properties implying \eqref{D_BdSv}.
The relation to the set of criteria given in \cite{AV_QSD}
and the applicability of our combined results in practice
is thus postponed to Subsection~\ref{D_sec_EC}.
Yet, it is also convenient to introduce or recall
concise and precise notations at use,
which we do first in Subsection~\ref{sec_nota}.
Subsection~\ref{sec_cpl_app}
is then devoted to the main contribution of this article,
namely the derivation of Property \eqref{D_BdSv}
thanks to our new criterion.
Subsection~\ref{D_sec_EC} is thus more generally concerned 
on the implications 
for the proof of exponential convergence to a unique QSD.
The above-mentioned complements to this result 
are given in Subsection~\ref{D_sec_Compl},
that are reciprocal statements  
and the uniformity in a localization procedure.

We detail our proofs in Section \ref{D_sec_Pf}.
Sections \ref{D_sec_Adapt3} and \ref{D_sec_Ad5} 
are finally dedicated to each of our applications,
which are to be introduced in their simpler form
in the next paragraph.

\subsection{Applications}
\label{D_sec_iApp}

The two illustrations of the current paper
are meant to help the reader 
get insight on the  adaptability
of our new criterion of ``almost perfect harvest".
They fall into the class of
piecewise deterministic Markov processes,
the second one being actually a pure jump process.
For the broader view of the applications 
we have in mind, 
we refer to the Discussion Section,
and more precisely Subsection \ref{sec_gen}.

We first consider 
the following process on $\bR^d$ for $d\ge 1$:
$$X_t 
= x - v\, t \mathbf{e_1}
+ \Tsum{i\le N_t} W_i,$$
where $x\in \bR^d$ is the deterministic initial condition,
$v>0$, $\mathbf{e_1}$ the first unit vector,
$N_t$ a Poisson process with intensity 1 
and the family $(W_i)_{i\in \bZ_+}$
is made up of i.i.d. normal variables 
with mean 0 and covariance matrix $\sigma^2 I_d$.
This process dies at a state-dependant rate.
The extinction event, whose time is denoted $\ext$,
is occurring at rate $t\mapsto \rho(X_t)$,
where $\rho(x):= \|x\|^2$ (with the euclidian norm).

Then, as stated in Theorem \ref{D_prop.A3},
there exists a unique quasi-stationary distribution associated 
to this dynamics and it is the only attractor of the conditioned marginals
as $t$ goes to infinity.

While the one-dimensional case is treated in \cite{CG20}, with a stronger result than in \cite{CH18},
the convergence result is new for the multidimensional setting.
For more details on the interpretation of this model 
and some generalization of its parameters, 
we refer to Section~\ref{D_sec_Adapt3}.
\\

The second illustration,
presented in Section \ref{D_sec_Ad5},
concerns a pure jump process on $\bR^d$, 
for $d\ge 2$:
\begin{equation*}
	X_t 
	:= x +  \Tsum{i\le N_t} \sigma W_i\cdot   \textbf{e}_{D_i}.
\end{equation*}
In this formula,
$x$ is the initial condition,
$N_t$ a standard Poisson process on $\bZ_+$,
while, for any $i\ge 1$,
$W_i$ is a standard 1d. normal random variable
and $D_i$ is uniform over $\II{1, d}$. 
Moreover, all these random variables are independent from each others.
Similarly as the previous example,
this process dies at a state-dependent rate 
given by $\rho_e : x\in \bR^d\mapsto \Ninf{x}^2$, 
where $\Ninf{x}:= \Tsup{i\le d} |x_i|$. 

Then, as stated in Theorem \ref{D_prop.Ad5},
$\sigma \le 1/8$ is a sufficient condition 
for the existence of a unique quasi-stationary distribution associated 
to this dynamics.

Note that jumps are restricted to happen 
along the vectors of an orthonormal basis 
$(\mathbf{e}_1,..., \mathbf{e}_d)$.
One jump is thus unsufficient to erase the singularities  
with respect to the Lebesgue measure on $\bR_d$.
To our knowledge, there is no other result of quasi-stationarity
for processes with such restrictions on the jumps. 
The result presented in \cite{CLW17} for the pure-jump case 
appears the closest to our.
In this multidimensional setting however, the jump effect is assumed to have a density
with respect to the Lebesgue measure on $\bR^d$.

More details on the interpretation of this model 
and some generalization of its parameters
are presented in Section \ref{D_sec_Ad5}.

\subsection{Specification of the results and techniques derived from \cite{AV_QSD}
and the current contribution}
\label{sec_spec}

\subsubsection{Current statement of the convergence we aim at}

For any Markov process $X = (X_t)_{t\in \bR_+}$ with extinction time $\ext$, 
the property of quasi-stationary convergence that we aim to prove 
for these processes is stated in the current article as follows.
The statement depends on a bounded function $\heig$ on the state space $\cX$
and on a probability measure $\alpha$ on $\cX$
 that are uniquely defined 
by imposing that $\int_{\cX} \heig(x) \alpha(dx)= 1$
in addition to the property.

There exists $C, \gamma>0$
such that the following inequality holds
for any $t\ge 0$
and  any probability measure $\mu$ on $\cX$
such that $\int_\cX \heig(x)  \mu(dx)>0:$
\begin{equation}
	\textstyle
	\|\, \PR_\mu \lp\, X_{\tp} \in dx;\; \tp < \ext \rp  - \alpha(dx) \,\int_\cX \heig(y)  \mu(dy) \|_{TV}
	\le C \|\mu - \alpha\|_{TV} e^{-\gamma \; \tp}.
	\label{intro_CVal}
\end{equation}
The reader is refered to Definition \ref{UEQS} in Section \ref{D_sec_ECV}
for a more precise statement
and to Subsection~\ref{sec_comp_QSD}
for the comparison with the result given in \cite{AV_QSD}.

\subsubsection{Main properties leading to the result following \cite{AV_QSD}}
Up to minor adjustment justified in Subsection~\ref{D_sec_PrCV}, 
it is the purpose of \cite{AV_QSD}
that the proof of this property \eqref{intro_CVal}
can be deduced 
by the proofs of the following 4 properties.
They depend on a probability measure $\alc$ on $\cX$
and on a subset $E$ of $\cX$,
that have to be the same.
The names are altered to ease the comparaison
with the assumptions proposed in the current article.
\begin{enumerate}
	\item[\textup{$(\overline{A0})$}] ``\textbf{Exhaustion of $\cX$}":
	There exists a sequence $(\cD_\ell)_{\ell\ge 1}$
	 of closed subsets of $\cX$ such that for any $\ell\ge 1$,
	  $\cD_\ell$ is included in the interior of $\cD_{\ell+1}$
	and such that their union makes up the whole state space $\cX$.
	\item[\textup{$(A1)$}] ``\textbf{Mixing property}":
For any $\ell\ge 1$, 
there exists $L>\ell$ and  $c, t>0$ such that:
\begin{equation*}
\frl{x \in \cD_{\ell}}
\hspace{.5cm}
\PR_x \lc {X}_{\tp}\in dx \pv
\tp < \ext \wedge T_{\cD_L} \rc 
\ge \cp\; \alc(dx),
\end{equation*}
where $T_{\cD_L}$ is the first exit time of $\cD_L$.

	\item[\textup{$(A2)$}] ``\textbf{Escape from the Transitory domain}": 
\textsl{There exists $\rho >\rho\iSv$ 
	such that the following boundedness property holds,
	where $\tau_E$ is the exit time of $E$
	and $E$ is required to be included in $\cD_\ell$ for $\ell$ sufficiently large:}
\begin{equation*}
\Tsup{x\in \cX} \;\E_{x} \lp 
	\exp\lc\rho\, (\ext\wedge \UDc) \rc \rp < \infty,
\end{equation*}
while the value $\rho\iSv$ is defined as a survival estimate with the following definition:
	\begin{equation*}
	\rho\iSv
	:= \sup\Big\{\gamma \in \bR\Bv 
	\sup_{L\ge 1} \inf_{t>0} \;
	e^{\gamma t}\,\PR_\alc(t < \ext\wedge T_{\cD_L}) 
	= 0
	\Big\}.
\end{equation*}

	\item[\textup{$(A3)$}] ``\textbf{Asymptotic Comparison of Survival}": 
	\begin{equation*}
	\limSInf{t}\; \sup_{x\in E} \; \dfrac{\PR_x(t<\ext)}{\PR_{\alc}(t<\ext)} 
	< \infty.
\end{equation*}
\end{enumerate}

\subsubsection{Current proposal of an alternative the last property}
\label{sec_prop}
Proving the three first properties 
is quite straightforward for the two presented applications,
as one can check in Subsections \ref{D_sec_th.A3}
and \ref{D_sec_th5}.
The proof of Property $(A3)$ is on the other hand much less direct
as it involves an asymptotic in large time for the proposed ratio.
It is only slighly easier than Equation \eqref{D_BdSv}
in that the inequality can be stated on a convenient subset of $\cX$
provided the exponential moment given in $(A2)$.

As mentioned just after Equation \eqref{D_BdSv} in the introduction,
neither the two-sided estimates
nor the gradient estimate 
or the Harnack inequality
are manifestly applicable.
This is due to the fact that any singularity of the initial condition
is maintained with just a decay in time
(and possibly a translation in space),
due to the jump event taking time.

This is why we make the emphasis on the following property 
that is quite easily proved in our examples.
It still depends on $E$ and $\alc$
but also on a value $\rho$ such that $\textup{$(A2)$}$ holds true
and on a value $\eps>0$
that is deduced from $(A1)$ and $(A2)$ 
as an intricate quantity.
In practice, 
we thus expect the property to be deduced whatever this value of $\eps$,
by adjusting the other estimates accordingly.
Notably, the property involves the design of a specific stopping time $\Uza$,
that we require to be infinite after a given threshold $\tZa$ in time. 
The design of the other stopping time $\UCa$ 
is to be adjusted according to $\Uza$, without specific restrictions.
\\

\textup{$(A3_F)$} ``\textbf{Almost Perfect Harvest}": 
\textsl{There exist $\tZa, \cp >0$ 
	such that for any $x \in E$
	there exists two stopping time 
	$\Uza$ and $\UCa$
	with the following properties:
\begin{equation*}
\PR_{x} \big(X(\Uza) \in dx' \pv \Uza < \ext \big) 
	\le \cp \,\PR_{\alc} \big(X(\UCa) \in dx'
	\pv \UCa < \ext\big),
\end{equation*} 
including the next conditions on $\Uza$:
\begin{equation*}
\PR_{x} (\Uza = \infty, \,  \tZa< \ext) 
	\le \fl\, \exp(-\rho\, \tZa),\qquad
	\text{ where }
	 \Lbr\ext \wedge \tZa \le \Uza \Rbr
	= \Lbr\Uza = \infty\Rbr.
\end{equation*} 
}

Additional regularity condition of $\Uza$
are also required with respect to the Markov property,
that we provide in detail in Subsection \ref{D_sec_AF3}.
Although technical, 
these conditions are expected to hold for any reasonable choice of $\Uza$,
generally defined through the meeting of specific conditions at some stopping time.

Thanks to Theorem \ref{D_Th:AF} in Subsection \ref{sec_main_thm},
assuming \textup{$(\overline{A0})$},
we deduce from $(A1)$, $(A2)$  and $(A3_F)$ that $(A3)$ holds also.

\subsubsection{Implications of the convergence results}
Thanks to Theorem \ref{D_AllPho},
which mostly relies on Theorems 2.1 and 2.2 of \cite{AV_QSD}
together with previous Theorem \ref{D_Th:AF},
the property of quasi-stationary convergence
is then deduced from $(A1)$, $(A2)$  and $(A3_F)$
(as well as from $(A1)$, $(A2)$  and $(A3)$).
We then infer in Corollary \ref{D_CVAl} the following convergence result for $\alpha$,
that justifies its identification as a quasi-stationary distribution:

``\textbf{Convergence to $\alpha$}":\textsl{ For any $t\ge 0$
	and  $\mu \in \M_1(\cX)$
	such that $\LAg \mu\bv \heig \RAg>0$:}
\begin{equation}
\|\, \PR_\mu \lc\, X_{\tp} \in dx \; 
| \; \tp < \ext \rc  - \alpha(dx) \, \|_{TV}
\le C \dfrac{\|\mu - \alpha\|_{TV}}
{\LAg \mu\bv \heig \RAg} \; e^{-\gamma \; \tp}.
\label{eq_CVAL}
\end{equation}	

It is actually justified in Theorem \ref{D_Th:AF} and \ref{D_AllPho}
that the family of sets $(\cD_\ell)$
may not cover the whole set $\cX$
for the new formulation given in Equation \eqref{intro_CVal} to hold.
Unexpectedly at first, this flexibility has been a great help 
for the proofs exploited in \cite{AV_GS} and \cite{AV_M},
even if $h$ is still strictly positive in the latter case.
In addition to the convergence depending on $\LAg \mu\bv \heig \RAg$,
it raises the question of identifying lower-bounds of $\heig$,
which we tackle in Proposition \ref{p:H0} by proving the following property:

``\textbf{Lower-bounds of $\heig$}":\textsl{ $\heig$ 
is uniformly bounded away from zero on any set $H\subset  \cX$
for which there exists $t>0$ and $\ell\ge 1$ such that 
$\quad 
\Tinf{x\in H}  
\PR_x(\tau_{\cD_\ell} < t \wedge \ext) > 0.$}

This justifies the identification of the domain of $\heig$ as follows:
$$\cH:= \{x\in \cX \pv h(x)>0\} = \Lbr x \in \cX\pv \Ex{\ell\ge 1}\PR_x(\tau_{\cD_\ell} < \ext)>0\Rbr.$$ 

As in \cite[Theorem 2.3]{AV_QSD},
we additionally deduce the existence of the $Q$-process,
that is of a Markov process living on this space $\cH$
whose generator $(\Q_x)_{x \in \cH}$
satisfies the following asymptotic property:
\begin{equation*}
	\limInf{\tp }
	\PR_x(\Lambda_\spr \bv \tp  < \ext) = \Q_x(\Lambda_\spr),
\end{equation*}
for any $x \in \cH$, $s>0$ and 
$\Lambda_\spr$ any $\F_\spr$-measurable set. 

As an extension of property \eqref{intro_CVal},
we also infer in Corollary \ref{D_QECV}
the following convergence property towards a unique stationary distribution $\beta$
of this $Q$-process.

``\textbf{Convergence to $\beta$}":\textsl{ 
	for any probability measure $\mu$ on $\cH$
	satisfying $\int_\cH \frac{\mu(dx)}{\heig(x)} < \infty$ 
	and $\tp \ge 0$:}
$$
 \NTV{ \Q_{\mu} \lc\, X_{\tp } \in dx\rc 
	- \beta(dx) }
\le C \,\; e^{-\gamma \; \tp} 
\cdot \left \|\mu(dx) - \left(\int_\cH \frac{\mu(dy)}{\heig(y)}\right)\cdot \beta(dx) \right\|_{1/\heig},
$$
where
$$\Q_\mu (dw):= \textstyle{\int_\cH}\mu(dx) \, \Q_x (dw),
\quad \|\mu\|_{1/\heig}:= \left \|\dfrac{\mu(dx)}{\heig(x)}\right \|_{TV}.$$

\subsubsection{Additional robustness properties of these results}

\paragraph{Approximation over restricted state space.}
The constants involved in the convergences are explicitly 
related to the parameters 
involved in the presented assumptions. 
Although the specific relation is very intricate,
it implies that one can fairly approximate the quasi-stationary regime
by restricting the state space $\cX$.
Indeed, for any $L\ge 1$ (thought to be large), 
let us consider the following approximation $\ext^L$ of the extinction event 
that restricts the process $X$ to remain in $\cD_L$
in the time-interval $[0, \ext^L]$.

Then, as stated in Theorem \ref{D_Approx}, 
by proving any of the two sets of assumptions (for the extinction time $\ext$) 
we deduce that all the above results hold uniformly in $L$
with the extinction time replaced by $\ext^L$.
We mean that we can choose the constants $C$ independent of $L$ 
for the convergences to $\alpha_L$, $\heig_L$, $(Q^L_t)$ and $\beta_L$.
Moreover, as $L$ goes to infinity, 
$\lambda_L$ converges to $\lambda$  and  $\alpha^L,h^L$ 
converge to $\alpha,h$  in total variation and pointwise respectively. 
Also, we deduce $\rho\iSv = \lambda$.

\paragraph{Reciprocal results.}
It is always satisfying to check that 
the sufficient conditions that one is attempting to verify 
are not too restrictive on the process.
This is why we have also been concerned with proving
of properties as analogous as possible to our key assumptions,
starting from the core result \eqref{intro_CVal}.
In particular, given this property of quasi-stationary convergence,
Proposition~\ref{D_recA2}
namely ensures that Property $(A2)$
holds for a certain parameter $\rho>\lambda$,
and a certain $\ell\ge 1$
 regardless of the choice of the sequence $\cD_\ell$ 
 satisfying Assumption \textup{$(\overline{A0})$}.
Proposition~\ref{D_recA3}
then asserts, 
under the same conditions and
for the same value of $\rho$,
that Property $(A3_F)$
is effectively met,
where $\zeta$ can be chosen as any probability measure on $\cX$.

\subsection{Comparison with the literature}
\label{D_sec_lit}

\subsubsection{Similar recent contributions}
\label{sec_comp_QSD}
The interest of our current result is mostly 
on generalizations of the Harris recurrence principle,
among  which 
\cite{ChQSD, CVly2, AV_QSD, BCGM19, CG20} present the most general statements
for homogeneous-in-time processes.
The upper-bounds that are derived 
	can be compared with the following inequality,
that takes a more general form than Equation~\ref{eq_CVAL}
	\begin{equation}\label{Cmu}
				\|\, \PR_\mu \left[\, X_{\tp} \in dx \; 
		| \; \tp < \ext \right]  - \alpha(dx) \, \|_{TV}
		\le C(\mu) \; e^{-\gamma \; \tp},
	\end{equation}
	where the main differences in the conclusions come from 
	different expressions and interpretations of the constant $C(\mu)$.
	\cite{ChQSD} first highlighted two necessary and sufficient conditions
	 for this convergence result to be uniform, that is with $C(\mu)$ taken as a constant.
The argument of a contraction in total variation norm is then simpler,
yet inspired the other extensions considering a non-trivial dependency $C(\mu)$
related to the forms of the key properties required.

When considering extinction, 
the property of linearity over the initial condition
is lost, so that linear expressions of $C(\mu)$ in terms of $\mu$
is not as natural as in the conservative case. 
This explains also why we came back to linear convergence statements 
in Equation~\ref{intro_CVal}.

The alternative techniques proposed in \cite{BCGM19}, \cite{CVly2} or \cite{CG20} involve 
contraction estimates of the operator 
in specific norms that are weighted by specific functions.
These functions satisfy some properties with respects to the semi-group $(P_t)$
that justify their association with the principles of Lyapunov contraction
in Harris recurrence techniques.
With the help of such Lypaunov function $W$ from $\cX$ to $\bR_+$,
we may generalize exponential moments as our property $(A2)$
and interpret the sequence $\cD_\ell$ as level sets of $W$ 
(i.e. $\cD_\ell=\{x\in \cX; W(x) \le \ell\}$).
With intricate boundary conditions, 
providing an efficient definition of such functions
in practice remains however challenging.

\subsubsection{Generalization of our approach}
\label{D_sec_gen}

Given the close interplay 
between our assumptions $(A2)$ and $(A3_F)$,
adapting the reasoning around $(A3_F)$ is not obvious.
Notably, the contraction estimates 
exploited in \cite{BCGM19} and \cite{CG20}
do not appear to relate as easily to a crucial property for our argument
(namely \eqref{D_eqSb} about the decay estimate of survival probability).

Doob's "$h$-transform" is a typical technique 
to deduce asymptotic results of a generally non-conservative semi-group $(P_t)_{t}$
from the study of a related sub-Markovian semi-group.
Provided that there exists a positive measurable function $\psi$ on $\cX$ and $\rho_\vee>0$ such that for any $t>0$, $P_t \psi \le e^{\rho_\vee t} \psi$,
the following definition 
indeed provides a sub-Markovian semi-group:
\begin{equation*}
	P^{\psi}_t(x, dy)
	:= \dfrac{\psi(y)}{\psi(x)} e^{\rho_\vee t} P_t(x, dy).
\end{equation*}
Implicitly, it means that $\psi$ is exploited to weight the state space $\cX$, as in the norms weighted by Lyapunov functions in \cite{BCGM19, CVly2, CG20}.
Note that the uniqueness property of the QSD 
given as 
$\psi(y) \alpha(dy)$
for $P^{\psi}$
corresponds exactly to the fact that $\alpha$ is the unique QSD of $P$
in the space $\M(\psi):= \{\mu \in \M(\cX)
\pv \LAg \mu\bv \psi\RAg < \infty\}$.
This allows for other probability distributions $\nu$ to be QSD,
in which case $\LAg \nu\bv \psi\RAg = \infty$ holds necessarily.

On the other hand,
our proofs would be easy to adapt to processes in discrete-time.
Our techniques should generalize naturally 
to time-inhomogeneous processes, 
given the recent adaptations 
presented in \cite{InhomChp, BenGenD, DV16, CG20}.
It can probably be extended to semi-Markov processes,
i.e. pure jump processes 
for which the waiting time between jumps 
is not necessarily exponential.

\subsubsection{Other frameworks}

General surveys 
like \cite{coll}, \cite{DP13} 
or more specifically 
for population dynamics \cite{MV12}
give an overview on the huge literature 
dedicated to QSD, 
for which Pollett has collected 
quite an impressive bibliography, 
cf. \cite{QSDbibli}. 

When jumps in continuous space are involved,
the reversibility property 
is generally expected not to hold true.
They are even more exceptional 
when the state space is multidimensional:
cf Appendix A of  \cite{CCM17}.

Comparison of survival is 
also an essential part of perturbation techniques as in Chapter 12 of \cite{Dm13} or in \cite{DM02},
yet it is mostly exploited for finite time
and compared to an intrinsic convergence rate. 
In \cite{FRS19},
results on the non-conservative semi-group
are deduced from the study of the Q-process 
as in the R-theory
(cf \cite{ANT80} for pure jump processes,
\cite{T74} or \cite{T74b} for the original discrete-time setting).
Our approach may provide guidance 
in dealing with estimates of the poorly known survival capacity. 

The other methods appear to bring less quantitative insight
in terms of uniqueness (except possibly \cite{ANT80})
or rate of convergence.
Besides the classical use of the Krein-Rutman theorem
(we recall \cite{CH18, HNV94}),
extensions from fixed point argument \cite{cmms}
and the above-mentioned R-theory, 
we also refer to ``renewal theory"
\cite{FKMP95}.

The compactness of the semi-group
is actually not required for our approach.
We recall that many classical approaches 
rely on this property to deduce the existence of a QSD
(cf e.g. the reviews \cite{coll}, \cite{MV12}),
often thanks to the Ascoli-Arzela theorem.
Since the process is allowed  in the illustrations given in Section \ref{D_sec_iApp}
either not to jump or to have a large number of jumps in any time-interval $[0, t]$,
we could not rely on this technique.

\section{Detailed description of our results}
\label{D_sec_MR}	
Before we present our results in more details, 
it is convenient to use efficient notations,
that we introduce or recall in the next subsection.
We can then focus in Subsection~\ref{sec_cpl_app} on the main contribution of the current article
before clarifying in Susbection~\ref{D_sec_EC} the implications of this result
combined with the ones of \cite[Theorems 2.1-3]{AV_QSD} on the quasi-stationary convergence.

\subsection{Notations}
\label{sec_nota}

In Subsection~\ref{D_sec_not} we describe our general notations, 
in Subsection~\ref{D_sec_abs} our specific setup
of a càd-làg strong Markov process with extinction
and clarify in Subsecton~\ref{sec_cond}
some notations made to express various event restrictions.

\subsubsection{Elementary notations}
\label{D_sec_not}

The most classical sets of integers are denoted as follows:
$\quad \bZ_+:= \Lbr 0,1,2...\Rbr$,$\; \bN:= \Lbr 1,2, 3...\Rbr$,
$\; [\![m, n ]\!]:= \Lbr m,\, m+1, ..., n-1,\, n\Rbr$ (for $m\le n$).
By the notation $:=$, we simply makes explicit that 
the equality is meant to explicit
some notation.
For maxima and minima, we use the following abbreviations:
$s \vee t:= \max\{s, t\}$,\,  
$s \wedge t:= \min\{s, t\}.$
In the paper, 
we may write $k\ge 1$ instead of $k\in \bN$
and $\tp \ge 0$ (resp. $c>0$)
instead $\tp \in \bR_+:= [0, \infty)$ 
(resp. $c\in \bR_+^* $ $:= (0, \infty)$),
when there is no real ambiguity.


The state space is denoted $\cX\cup {\partial}$,
where the cemetery $\partial$ is assumed to be isolated from the topology $\B$ 
on the Polish space $\cX$.
In the study of the process,
we will need to apply the Markov property at first entry times 
of either closed or open subsets.
This is why we assume that the time homogeneous process $X$
is strong Markov
for the filtration $(\F_\tp)_{\tp\ge0}$ that is right-continuous and complete
and that it has càd-làg paths (left limited and right continuous).
The hitting time (resp. the exit time out) of $\cD$,
for some domain $\cD \subset \cX$, 
will generally be denoted by $\tau_{\cD}$ (resp. by $T_{\cD}$).
These are stopping times for any $\cD$ that is either closed or open, 
cf. Theorem~52 in \textup{\cite{Mey}},
or more recently Theorem~2.4 in \textup{\cite{B10}}.

Exploiting the same notations as in \cite{Rog00}, Definition III.1.1,
$P_\tp$ would then be the semi-group of the process 
and the latter shall satisfy the usual measurability assumptions 
and the Chapman-Kolmogorov equation.
The law of the process starting from initial condition $x\in \cX\cup {\partial}$
will be given by the probability $\PR_x$.

%

\subsubsection{The stochastic process with absorption}
\label{D_sec_abs}

Here, we consider a strong
Markov processes absorbed at $\partial$: the cemetery. 
More precisely, we assume that 
$X_s = \partial$ implies $X_\tp = \partial$ for all $\tp \ge s$. 
This implies that the extinction time:
$\quad \ext:= \inf\Lbr \tp \ge 0 \pv X_\tp = \partial \Rbr\quad$
is a stopping time. 
Thus, we rather consider the family $(P_\tp)_{\tp \ge0}$
as a non-conservative semigroup of operators on the set \mbox{$\B_+(\cX)$}
(resp.  \mbox{$\B_b(\cX)$}) 
of positive (resp. bounded)
$(\cX,\B)$-measurable real-valued functions. 

For any probability measure $\mu$ on $\cX$, 
and \mbox{$f\in \B_+(\cX)$} 
(or \mbox{$f\in \B_b(\cX)$}), 
we use the notations:\\
$$\PR_\mu (.):= \int_{\cX} P_x(.) \; \mu(dx), \quad
\langle \mu \bv f\rangle:= \int_{\cX} f(x) \; \mu(dx).$$
We also denote by $\E_x$ (resp. $\E_\mu$) 
the expectation according to $\PR_x$ (resp. $\PR_\mu$).

The set of respectively probability measures on $\cX$, of positive and of signed measures
are denoted respectively $\Mone$, $\M_+(\cX)$ and $\M(\cX)$.
For any $B\in \B$ and $\mu \in \Mone$, 
$\mu P_t(B)$ is clearly defined as $\PR_\mu(X_t \in B)$.
The fact that $t<\ext$ immediately follows from $B$ being a subset of $\cX$
and $X$ being absorbed in $\partial$.
Yet, we wish to avoid confusion for instance in the examples 
given in Subsection~\ref{D_sec_iApp},
where extinction is prescribed through some state-dependent rate 
defined a posteriori 
from an internal dynamics of $X$ on $\cX$.
This is why we usually make the restriction on the event $t<\ext$ explicit 
even when it is not required, notably 
in the definition of the action of $P_t$ on $\mu\in \Mone$ as follows:
\begin{equation*}
\mu P_\tp(dy):= \PR_\mu(X_\tp \in dy\pv t<\ext),
\quad \text{ or }
\langle \mu P_\tp\bv f\rangle 
\,=\, \langle \mu \bv P_\tp f\rangle 
\,= \E_\mu[f(X_\tp)\pv t<\ext],
\end{equation*}
where $f\in \B_+(\cX)$ or \mbox{$f\in \B_b(\cX)$}.
Let us then define the family of conditioned operators $(A_t)_{t\ge 0}$
acting as follows on any probability measure $\mu \in \Mone$:
\begin{equation*}
	\mu A_\tp(dy):= \PR_\mu(X_\tp \in dy \bv \tp < \ext), \qquad 
\langle\mu A_\tp\bv f\rangle \,= \E_\mu[f(X_\tp)\bv \tp< \ext]
	= \dfrac{\E_\mu[f(X_\tp)\pv t<\ext]}{\PR_\mu[t<\ext]}.
\end{equation*}
$\mu A_\tp$ is what we call the MCNE
at time $\tp$, with initial distribution $\mu$,
as it is "the Marginal distribution (at time $t$)
Conditioned upon the fact 
that No Extinction has yet occurred"
(also at time $t$). 

In this setting, 
the family $(P_\tp)_{\tp\ge0}$ 
(resp. $(A_\tp)_{\tp\ge0}$) 
defines a linear but non-conservative semigroup
(resp. a conservative but non-linear semigroup)
of operators on $\Mone$ 
endowed with the total variation norm. 
We consider the following definition of this norm,
generally for any signed measure $\mu \in \M(\cX)$:
$ \quad
\|\mu \|_{TV}:= \sup\Lbr |\mu(A)|\pv A \in \B \Rbr.
\quad $
While the semi-group $P$ is directly generalized by linearity for any signed measure,
note that it is not as clear for the semi-group $A$ because $\PR_\mu[t<\ext]$ could then be equal to zero for some $t>0$.

A probability measure $\alpha$ is said to be the \textit{quasi-limiting distribution} of an initial condition $\mu\in \Mone$ if:
$\quad 
\frl{B \in \B}\quad
\limInf{\tp} \PR_\mu(X_\tp \in B\bv \tp< \ext):= \limInf{\tp} \mu A_\tp(B) = \alpha(B).
\quad $
\\
\noindent
It is now classical 
(cf e.g. Proposition~1 in \cite{MV12}) 
that $\alpha$ is then a quasi-stationary distribution or QSD, in the sense that:
$\frl{\tp\ge 0}\quad
\alpha A_\tp(dy)
= \alpha(dy).$
Also, for any such QSD, 
there exists a unique extinction rate $\laZ>0$ 
such that: $\frl{t \ge 0} \PR_\alpha(t< \ext) = \exp[-\laZ\, t]$.


\subsubsection{Conditions, stopping times and random events}
\label{sec_cond}
While dealing with the Markov property between different stopping times, 
we wish to clearly indicate with our notation 
that we introduce a copy of $X$ 
(ie a process with the same semigroup $P_\tp$)
independent of $X$ given its initial condition. 
This copy 
(and the associated stopping times) 
is then denoted with a tilde 
($\wtd{X},\, \wtd{\ext},\, \wtd{T}_\cD$ etc.). 
For instance in the notation 
$\PR_{X(\UDc)} (t- \UDc < \wtd{\ext})$, 
$\UDc$ and $X(\UDc)$ refer to the initial process $X$ 
while $\wtd{\ext}$ refers to the \mbox{copy $\wtd{X}$.}

Besides, some notations of semi-colons  and commas 
have meanings that are specific to our probabilistic notations
and currently very efficient given that we often consider restrictions on various events.
For expectations,
the terms after the semi-colon indicate the sets of conditions
under which the former term is counted 
(it is replaced by 0 if the conditions are not met).
For sets, and notably random ones,
the terms after the semi-colon indicate the sets of conditions
under which the former term is included.
Different conditions may be separated by commas,
notably when combinations of those are introduced through indices.
For instance, 
given any random variables $F$, $(X_i)_{i\le d}$, for $d\ge 1$, $T, T'$ and $t>0$,
the following notation can be translated as follows:
\begin{equation*}
	\E(F\pv \frl{i\le d}X_i\ge 0\mVg t<T) 
	=\textstyle \E\big(F\cdot \prod_{\{i\le d\}} \idc{X_i\ge 0} \cdot \idc{t<T}\big).
\end{equation*}
To give another example, the following notation:
\begin{equation*}
	\{s\ge T'\pv  \frl{i\le d} X_i\ge 0\mVg s<T\}
\end{equation*}
is to be understood as the random set $E$
(that is thus dependent on $\om \in \Omega$)
defined as follows.
If there exists $i$ such that $X_i(\om)$ is negative or if $T'(\om)\ge T(\om)$, then $E =\emptyset$.
Otherwise, $E$ consists
of all values $s\in \bR_+$
such that $s\in [T'(\om), T(\om))$.
The infimum of an empty set is generally to be taken as $\infty$.
Without such semi-colon, the set has to be considered 
as the part of $\Omega$ for which the conditions are satisfied, 
for instance:
\begin{equation*}
	\{\frl{i\le d}X_i\ge 0\mVg t<T\}
	= \cap_{\{i\le d\}}\{X_i\ge 0\}\cap \{t< T\}
\end{equation*}
consists of all $\om \in \Omega$ such that 
for any $i\le d$, $X_i(\om)\ge 0$ and $t<T(\om)$.

\subsection{Coupling approach including failures}
\label{sec_cpl_app}

The main contribution of this article 
is presented in this Subsection~\ref{sec_cpl_app},
namely the derivation of Property \eqref{D_BdSv}
thanks to our new criterion.
We start in Subsection~\ref{sec_key_ass} 
by explaining the key properties 
derived from \cite[Section 2]{AV_QSD} and the role of the associated parameters
involved in our new criterion.
After presenting this new criterion of "Almost Perfect Harvest" in Subsection~\ref{D_sec_AF3},
we state our main Theorem \ref{D_Th:AF} in the next Subsection~\ref{sec_main_thm}.

\subsubsection{The associated assumptions}
\label{sec_key_ass}

Our core property of "almost perfect harvest"
exploits several parameters ($\rho > 0$, 
$E$ and $\alc$)
whose choices are strongly tied to 
several key properties highlighted in \cite{AV_QSD}.
We shall thus start 
by explaining the key properties of \cite{AV_QSD} on which we rely
and thus the role of these parameters.
\\

The approach is trajectory-based 
and designed to handle specific dependencies on the initial condition,
so it appears efficient
to consider a customizable covering 
by an increasing sequence of sets, as specified in the next property.
\\

$(\underline{A0})$ : \textbf{``Specification of space"}
\quad
There exists a sequence $(\cD_\ell)_{\ell\ge 1}$ of closed subsets of $\cX$ such that for any $\ell\ge 1$, $\quad  \cD_\ell \subset 
	int(\cD_{\ell+1}) \quad$
	(with $int(\cD)$ the interior of $\cD$).
	\\

	This sequence serves as a reference for the other key properties, 
	notably through the following notation,
	for subsets that are "regular" with respect to this specification:
	\begin{equation}
		\Dps:= \Lbr \cD\pv \cD \text{ is closed and } 
		\Ex{\ell\ge 1} \cD \subset \cD_\ell\Rbr.
		\label{D_Dps}
	\end{equation}
	
	As it is stated in  \cite[Section 2]{AV_QSD} (and by extension in the introduction),
	this property \textup{$(\underline{A0})$} is often strengthened as follows.
	This additional condition of complete covering
	enables to obtain uniqueness of the QSD as such
	and not only among the QSD with minimal extinction rate.
\\
	
	\textup{$(\overline{A0})$} : the sequence  $(\cD_\ell)_{\ell\ge 1}$ satisfying \textup{$(\underline{A0})$}
	is such that $\quad \medcup_{\ell\ge 1} \cD_\ell= \cX$.
	\\
	
The next assumption consists in a minoration of the density of the process
after a given time. This minoration extends the classical Doeblin inequality
where the marginal law is lower-bounded uniformly over the initial conditions.
This is the step that produces the mixing of the past dependencies.
As in Harris recurrence technique, 
this lower-bound needs only to be uniform locally in the initial condition,
at the expense of an associated contraction estimate.
In our case, we further restrict the probability to trajectories that remain locally confined
and exloit for practical convenience
the minoration with the restriction on the event of survival $\{t<\ext\}$
rather than with the conditioning on the event.
	We recall the following definitions 
	for the exit and first entry times of any set $\cD$:
	$$T_{\cD}:= \inf\Lbr  \tp \ge 0 \pv X_\tp \notin \cD \Rbr
	\mVg\quad
	\tau_\cD:= \inf\Lbr \tp \ge 0\pv X_\tp \in \cD \Rbr.$$
	
	\noindent
	\textup{$(A1)[\alc]$} : \textbf{``Mixing property"} \quad
	The probability measure $\alc \in  \Mone$ is such that, 
	for any $\ell\ge 1$, 
	there exists $L>\ell$ and  $c, t>0$ such that:
	\begin{equation*}
 \frl{x \in \cD_{\ell}}
		\hspace{.5cm}
		\PR_x \lc {X}_{\tp}\in dx \pv
		\tp < \ext \wedge T_{\cD_L} \rc 
		\ge \cp\; \alc(dx).	
	\end{equation*}
	
	The last assumption corresponds to the contraction estimate in Harris recurrence techniques,
	in the form of an exponential moment of return to some reference subset $E$ of $\cX$.
	To combine this property with the mixing stated in $(A1)$,
	we require this set $E$ to be included in one of the $\cD_\ell$.
	For practical convenience,
	we consider the moment estimate on the stopping time $\ext\wedge \UDc$
	for the unconditioned expectation.
	For the contraction to be of relevance for the conditioned process,
	an exponential moment estimate with a parameter larger 
	than the extinction rate of the process appears crucial
	the more so as the extinction rate is not directly accessible
	and has to be evaluated.
	\\
	
	\textup{$(A2)[\rho, E]$} : \textbf{``Escape from the Transitory domain"} \quad
	\textsl{The value $\rho >0$ and the subset $E \in \Dps$ are such that:}
	\begin{equation*}
\Tsup{x\in \cX} \;\E_{x} \lp 
		\exp\lc\rho\, (\ext\wedge \UDc) \rc \rp < \infty.
	\end{equation*}
	$\rho$ as stated in $(A2)$ and $(A3_F)$ 
	is required in the following theorems
	to be strictly larger
	than the following \textbf{``survival estimate"}: 
	\begin{equation*}
		\rho\iSv
		:= \sup\big\{\gamma \in \bR\Bv 
		\sup_{L\ge 1} \inf_{t>0} \;
		e^{\gamma t}\,\PR_\alc(t < \ext\wedge T_{\cD_L}) 
		= 0
		\big\}.
		\label{D_rsv}
	\end{equation*}
	Moreover, the set $E$ shall be common for $(A2)$ and $(A3_F)$.
	
	  	From the proof of Corollary 5.2.1
	in \cite{AV_QSD},
	we deduce the following property of survival
	as a direct consequence the above properties.
	
	\begin{prop}
Assume that \textup{$(\underline{A0})$}, \textup{$(A1)[\zeta]$} and \textup{$(A2)[\rho, E]$}
hold for some $\zeta\in \Mone$, $\rho>\rho\iSv$ and $E\in \Dps$. 
Then, there exists $\tp\iSB, c\iSB>0$ and $\Rsv \in (\rho\iSv, \rho)$ such that:
\begin{equation}
\frl{ u\ge 0}
\frl {\tp \ge u + \tp\iSB}
\hspace{0.5cm}
P_{\alc}(\tp - u < \ext) 
\le \cp\iSB\; e^{\Rsv\, u}\;  \PR_{\alc}\lp  \tp < \ext \rp.
\label{D_eqSb}
\end{equation}
	\end{prop}
It means that for sufficiently large times,
the exponential decay in time with parameter $\Rsv$ 
provides a relevant estimation 
of the decay of the probability of survival starting from the initial conditon $\zeta$.

	\subsubsection{The new key property}
	\label{D_sec_AF3}
	
	We are now in position to
	introduce the following weak form of Harnack inequality,
	whose purpose is to imply the property 
	\eqref{D_BdSv}:\\
	
	\textup{$(A3_F)[\rho, E, \alc]$}~: \textbf{``Almost perfect harvest"}\quad
	\textsl{The value $\rho >0$, the subset $E \in \Dps$ 
		and the probability measure $\alc \in \Mone$ are such that the following condition holds:\\
		For any $\fl\in (0,\, 1)$,
		there exist $\tZa, \cp >0$ 
		such that for any $x \in E$
		there exists two stopping time 
		$\Uza$ and $\UCa$
		with the following properties:} 
\begin{equation}
\PR_{x} \big(X(\Uza) \in dx' \pv \Uza < \ext \big) 
\le \cp \,\PR_{\alc} \big(X(\UCa) \in dx'
\pv \UCa < \ext\big),
\label{D_Abs}
\end{equation} 
including the next conditions on $\Uza$:
\begin{equation}
\Lbr\ext \wedge \tZa \le \Uza \Rbr
= \Lbr\Uza = \infty\Rbr
\quad  \text{ and }\quad
\PR_{x} (\Uza = \infty, \,  \tZa< \ext) 
\le \fl\, \exp(-\rho\, \tZa).
\label{D_FL}
\end{equation} 
	\textsl{As a regularity condition of $\Uza$ with respect to the Markov property,
		we additionally require that
		$\Uza$ can be expanded in the following sense
		by at least some stopping time $U_{H}^\infty$ such that:}
	
	\textsl{$\star$   $U_{H}^\infty:= \Uza$ 
		on the event  $\Lbr \ext\wedge \Uza < \tau\iET^1\Rbr$, 
		where $\tau\iET^1:= \inf\{s\ge \tZa \pv
		X_s \in E\}$. }
	
	\textsl{$\star$ On the event 
		$\Lbr \tau\iET^1 \le \ext \wedge  \Uza \Rbr$
		and conditionally on $\F_{\tau\iET^1}$,
		the law of $U_{H}^\infty - \tau\iET^1$ coincides
		with the one of $\wtd{U}_{H}^\infty$ 
		for a realization $\wtd{X}$ 
		of the Markov process $(X_t, t\ge 0)$
		with initial condition $\wtd{X}_0:= X(\tau\iET^1)$ 
		and independent of $X$ conditionally on $X(\tau\iET^1)$.}

\paragraph{Interpretation of the core assumption\\}

This property is to be compared to the more classical Harnack inequality,
which should take the following form (as can be seen in \cite{AV_QSD}):\\
\textsl{The subset $E \in \Dps$ 
	and the probability measure $\alc \in \Mone$ are such that 
	there exist $\tZa, v, \cp >0$ for which the following inequality holds 
for any $x \in E$:}
\begin{equation}
	\PR_{x} (X_{\tZa} \in dx' \pv \tZa < \ext ) 
	\le \cp \,\PR_{\alc} (X_v \in dx'
	\pv v < \ext).
	\label{har}
\end{equation}

Thanks to \textup{$(A3_F)$},
we will also be able to couple the trajectories of the processes starting with the initial conditions
resp. $x\in E$ and $\alc$ so that they coincide with a time lag for sufficiently large times.
In other words, we want to embed the ``tail trajectories" of the process starting from $x$
 into the trajectories starting from $\alc$, a time lag being allowed.
When \eqref{har} holds true, the trajectories with initial condition $x$
can namely be coupled after time $\tZa$, with a time-lag of $\tZa - v$.

 This embedding procedure is what we describe as the ``harvest"
 of the tail trajectories.
With $(A3_F)$, this embedding is here to occur after time $U_{H}^{\infty}$,
the ``harvesting time",
the time-lag being more flexible and notably allowed to be random.

The Markov property being granted, inequality \eqref{D_Abs} is what makes this embedding possible. 
This embedding is directly obtained
on the event $\{U_H<t_F\wedge \ext\}$,
in which case $U_H$ and $U_H^\infty$ coincide.
The Markov property being exploited for the coupling, 
we request both $U_{H}$ and $V$ to be stopping time.

It is however rarely expected that $U_H$ can be obtained upper-bounded by a uniform constant $t_F$
 such that inequality \eqref{D_Abs} still holds.
 Possibly several attempts of coupling 
 may thus be required under some specific conditions, 
 that we describe as failures.
A ``failure" to do the coupling
is then characterized by the event  $\{\ext \wedge \tZa \le \Uza\}$
 (which we require to coincide with $\{\Uza = \infty\}$).
The objective of establishing the embedding of tail trajectories 
except in case of rare events in probability
is what motivates our denomination of "almost perfect harvest".

In case of failures,
$U_{H}^{\infty}$ must be larger than $\tZa$,
and a new attempt can directly be implemented only after $\tau\iET^1$.
We need to ensure that the events of recurrent failures play a negligible role
in the probability of long-term survival.
This is why,
in complement to the factor $\epsilon$ being sufficiently small,
we require the penalisation by $\exp(-\rho\, \tZa)$
to compensate for the extinction rate 
(during the time-interval $[0, \tZa]$ for the process with initial condition $\zeta$).
For such a compensation to be exploitable, 
the value of $\rho$ will be required to be greater 
than the known lower bound $\rho_S$
 for the decay of the survival probability.

The condition of failure where $\{U_{H} = \infty\}$
has to be adjusted for convenience to the model. 
As the allowed error $\fl$ goes to 0,
we expect to see the failure conditions stated 
through certain thresholds 
being less and less stringent.
This certainly leads to larger and larger constants $c$ in \eqref{D_Abs},
and often to larger and larger constants $t_F$.

Because the attempts of coupling are meant to be iterated after each failure,
the condition related to the Markov property on $U_H^\infty$
is requested.
It takes into account the waiting time before  the process comes back to $E$
and a new attempt can be initiated.
Provided  that the definition of $U_H$ do not depend in a singular way 
on the initial condition $x\in E$,
this condition on $U_H^\infty$ should be easily satisfied.

%
%

\subsubsection{The central result of the paper}
\label{sec_main_thm}

\begin{theo}
	\label{D_Th:AF}
	Assume that there exist $\alc\in~\Mone$,
	$\tp\iSB, c\iSB>0,$
	$\rho \ge \Rsv > \rho\iSv >0$,  and  $E \in \Dps$ such that  
	inequality \eqref{D_eqSb},
	assumptions \textup{$(\underline{A0})$}, \textup{$(A1)[\alc]$}, \textup{$(A2)[\rho, E]$}
	and  \textup{$(A3_F)[\rho, E, \alc]$}  hold. 
	Then:
	\begin{equation}
\limSInf{t}\; \sup_{x\in E} \; \dfrac{\PR_x(t<\ext)}{\PR_{\alc}(t<\ext)} 
< \infty.
\label{D_A3}
	\end{equation}
\end{theo}

For the implication that 
\eqref{D_BdSv} is a consequence
($(\underline{A0})$, $(A1)$ and $(A2)$ being granted)
of this more local property
(which is exactly $(A3)$),
we refer to \cite[Theorem 5.2]{AV_QSD}.
A careful check of the proof ensures
that  the additional requirement $\medcup_{\ell\ge 1} \cD_\ell= \cX$
is actually not exploited therein.

\begin{rem}
	\label{RConf}
In $(A1)$,
a confinement in some subspace $\cD_L$
is required. 
This confinement part of the assumption
is actually not involved in the proof of Theorem \ref{D_Th:AF},
but is 
in fact exploited to deduce Theorem \ref{D_AllPho}.
We kept it to avoid too frequent variations
of our assumptions.
\end{rem}

\subsection{Exponential convergence to a unique $QSD$}
\label{D_sec_EC}

Thanks to Theorem \ref{D_Th:AF},
we see that the general conclusions of \cite[Theorems 2.1-3]{AV_QSD}
can be derived from a new set of assumptions,
as we present it in this Subsection.

 \subsubsection{Two general sets of assumptions}

We say that Assumption $\mathbf{(A_F)}$ holds, whenever:

``There exists $\alc\in \Mone$ such that $(A1)[\alc]$ holds 
for a specific sequence $(\cD_\ell)$ satisfying $(\underline{A0})$.
Moreover, there exists $\rho  > \rho\iSv$ 
and $E \in \Dps$ such that
assumptions $(A2)[\rho, E]$
and $(A3_F)[\rho, E, \alc]$  hold."
\\

Slightly adapting \cite{AV_QSD} (regarding $(\underline{A0})$), we say that Assumption $\mathbf{(A)}$ holds, whenever:

``There exists $\alc$ such that $(A1)[\alc]$ holds 
for a specific sequence $(\cD_\ell)$ satisfying $(\underline{A0})$.
Moreover, there exists $\rho  > \rho\iSv$ 
and $E \in \Dps$ such that
assumptions $(A2)[\rho, E]$
and $(A3)[E, \alc]$  hold."

\begin{rem}
	\label{rm.CV}
	$\star$ The purpose of Theorem \ref{D_Th:AF} is to prove 
	that Assumption $\mathbf{(A_F)}$ 
	actually implies Assumption $\mathbf{(A)}$.
	For this, we refer to Subsection \ref{D_sec_PrCV}.

	$\star$ Almost sure extinction is not at all needed for our proof (which would in fact include the case where there is no extinction, or only in some ``transitory domain").
	
	$\star$	In our set of assumptions, the case 
	where  $\cD_\ell:= \cX$ for any $\ell$
	is actually allowed 
	and fits in the setting of uniform convergence
	covered in \cite{ChQSD}.
\end{rem}

%
%
%
%
%
%
%
%
%

%

\subsubsection{Main consequences of our main result}
\label{D_sec_ECV}

The main property that we deduce in Theorem \ref{D_AllPho} 
and in our subsequent applications
can be stated as in the following definition:
\begin{defnt}
	\label{UEQS}
	For any linear positive and bounded semi-group $(P_t)_{t\ge 0}$
	acting on a state space $\cX$, 
	we say that $P$ displays a uniform exponential quasi-stationary convergence
	with characteristics
	$(\alpha, h, \lambda) \in \M_1(\cX)\ltm B(\cX)\ltm \bR$
	if $\LAg \alpha \bv h\RAg = 1$
	and there exists $C, \gamma>0$ such that for any $t>0$
	and for any measure $\mu\in \M(\cX)$ with $\NTV{\mu}\le 1$:
	\begin{equation}
	\NTV{e^{\lambda t} \mu P_t- \LAg \mu\bv h\RAg \alpha}
	\le C e^{-\gamma t}.
		\label{PtCV}
	\end{equation}
\end{defnt}
This definition slightly extends the one given in introduction, cf \eqref{intro_CVal},
by considering a general class of semi-groups
and allows the value $\lambda$ to be negative.
Since $\LAg \alpha \bv h\RAg = 1$ and since $P_t$ is linear,
elementary computations show that the above property
implies the apparently stronger one given in \eqref{intro_CVal}
with $\NTV{\mu - \alpha}$ as additional factor.

The property given in Definition \ref{UEQS} implies 
that $\alpha$ is a QSD with extinction rate $\lambda$,
as stated in the following fact, whose proof is deferred 
to the appendix.
\begin{fact}
	\label{fQSD}
If a semi-group $P$ displays a uniform exponential quasi-stationary convergence
with characteristics
$(\alpha, h, \lambda)$, 
then for any $t>0,$ $\alpha P_t(ds) = e^{-\lambda t} \alpha(ds)$.
\end{fact}

\begin{rem}
$\star$	It is elementary that
 $h_t:x\mapsto e^{\lambda t} \LAg \delta_x P_t\bv \idg{} \RAg$
  converges in the uniform norm to $h$,
  implying that $h$ is non-negative.
  Since $h_{t+t'} = e^{\lambda t} P_t h_{t'}$, 
	one can also easily deduce that $e^{\lambda t} P_t h = h$.
%
	
$\star$	 By ``characteristics", we express that they are uniquely defined.	
\end{rem}

\begin{theo}
\label{D_AllPho}
Assume that either $\mathbf{(A_F)}$ or $\mathbf{(A)}$ holds.
	Then, 
	$P$ displays a uniform exponential quasi-stationary convergence
	with some characteristics
	$(\alpha, h, \lambda)$ in $\M_1(\cX)\ltm B(\cX)\ltm \bR_+$.
%
%
It also holds true that $h$ is bounded away from zero on any $\cD_\ell$, 
that $h$ is upper-bounded on $\cX$
and that $\alpha(\cD_{\ell})>0$ for $\ell$ sufficiently large.
	 \end{theo}
 

The convergence to $\alpha$ is made more precise by the following corollary:
\begin{cor}
\label{D_CVAl}
Assume \eqref{PtCV}.
Then for any $t\ge 0$
and  $\mu \in \M_1(\cX)$
such that $\LAg \mu\bv \heig \RAg>0$:
	\begin{equation}
	\|\, \PR_\mu \lc\, X_{\tp} \in dx \; 
	| \; \tp < \ext \rc  - \alpha(dx) \, \|_{TV}
	\le C \dfrac{\|\mu - \alpha\|_{TV}}
	{\LAg \mu\bv \heig \RAg} \; e^{-\gamma \; \tp}.
	\label{D_ECvAl}
\end{equation}	
\end{cor}

Thanks to Theorem \ref{D_AllPho}, it is elementary 
	that $h$ is positive under assumption $(\overline{A0})$. 
	It might be useful otherwise to exploit the following proposition 
	to identify a posteriori the domain on which $h$ is positive.

\begin{prop}
	\label{p:H0}
	Assume that $\mathbf{(A_F)}$ or $\mathbf{(A)}$ holds.
Then, the survival capacity $h$ 
is uniformly bounded away from zero on any set $H\subset  \cX$
that satisfies the following condition:\\
$(H_0) :$	there exists $t>0$, $\ell\ge 1$ such that 
$\quad 
	\Tinf{x\in H}  
	\PR_x(\tau_{\cD_\ell} < t \wedge \ext) > 0.$
\\
	It implies the following identification 	:
	$$\cH:= \{x\in \cX \pv h(x)>0\} = \Lbr x \in \cX\pv \Ex{\ell\ge 1}\PR_x(\tau_{\cD_\ell} < \ext)>0\Rbr.$$
\end{prop}

\begin{rem}
Corolllary \ref{D_CVAl}
implies that if $\nu$ were a QSD different from $\alpha$,
then $\nu A_t = \nu\neq \alpha$ 
thus $\LAg \nu \bv h\RAg = 0$.
This implies
$\nu(\cH ) = 0$.
So by contraposition, 
 the previous proposition
 provides a practical way to ensure 
 the uniqueness of the QSD a posteriori.
\end{rem}

Once the set $\cH$ is clarified, 
we can study Doob's $h$-transform 
of the semi-group $P$ with weight given by the survival capacity.
This is the so-called $Q$-process that is actually a conservative Markov process as stated in the next corollary:

\begin{cor}
\label{D_QECV}
Under again $\mathbf{(A_F)}$ or $\mathbf{(A)}$,
 with $(\alpha, h, \lambda)$ 
the characteristics of exponential convergence of $P$, 
the following properties hold:

\noindent
 \textup{(i) \textbf{Existence of the $Q$-process:}} 
 
There exists a family $(\Q_x)_{x \in \cH}$ of probability
measures on $\Omega$ defined by:
\begin{equation}
\limInf{\tp }
\PR_x(\Lambda_\spr \bv \tp  < \ext) = \Q_x(\Lambda_\spr),
\label{Qxdef}
\end{equation}
for any $x \in \cH$, $s>0$ and 
$\Lambda_\spr$ any $\F_\spr$-measurable set. 
The process $(\Omega;(\F_\tp )_{\tp \ge 0};
(X_\tp )_{\tp \ge 0};(\Q_x)_{x \in \cH})$
is an $\cH$-valued homogeneous strong Markov process.
\end{cor}

\noindent
\textbf{(ii) Weighted exponential ergodicity of the $Q$-process:}
	
\textit{The measure $\, \beta(dx):= \heig(x)\, \alpha(dx)$ 
is the unique invariant probability measure under $\Q$.
\\
Moreover, for any $\mu \in \M_1(\cH)$
satisfying $\LAg \mu\bv 1/\heig\RAg < \infty$ 
and $\tp \ge 0$:}
 	\begin{equation}
 \NTV{ \Q_{\mu} \lc\, X_{\tp } \in dx\rc 
 - \beta(dx) }
 	\le C \,\|\mu - \LAg\mu\bv 1/h \RAg \, \beta\|_{1/\heig}\; e^{-\gamma \; \tp},
 		\label{D_ECvBeta}
 	\end{equation}
 where $\Q_\mu (dw):= \textstyle{\int_\cH}\mu(dx) \, \Q_x (dw),
 \quad \|\mu\|_{1/\heig}:= \left \|\dfrac{\mu(dx)}{\heig(x)}\right \|_{TV}.$

The constant $\LAg \mu\bv 1/h\RAg $ before $\beta$ \eqref{D_ECvBeta}
is optimal up to a factor 2,
due to the following fact,
that is proved in the appendix.
This is also the case for the implicit constant 1
before $\alpha$ in \eqref{D_ECvAl}.
 
\begin{fact}
	\label{D_RemOpt}
For any $u>0$:
$$\|\mu - u\,\alpha\|_{1/h} 
\ge  (1/2) \cdot\|\mu - \LAg\mu\bv 1/h \RAg \beta\|_{1/h} ,
\qquad \|\mu - u\,\alpha\|_{TV} 
\ge (1/2) \cdot\|\mu - \alpha\|_{TV}.$$
\end{fact}

\subsubsection{Additional robustness properties of the results}
\label{D_sec_Compl}

\paragraph{Reciprocal results}
The following propositions are stated to justify 
that our assumptions are not particularly restrictive.
We shall see in the following Theorem \ref{D_Approx}
that $\rho\iSv$ actually equals $\lambda$
if  $\mathbf{(A_F)}$ is satisfied.
We thus aim at stating reciprocal statements 
in terms of a parameter $\rho$ satisfying $\rho > \lambda$.
The first proposition concerns $(A2)$,
that is derived under the assumption of exponential convergence
together with property \textup{$(\overline{A0})$}.
\begin{prop}
	\label{D_recA2}
	Assume that $\cX = \medcup_{\ell\ge 1} \cD_\ell$ where $(\cD_\ell)_{\ell\ge 1}$ is an increasing subsequence.
	Assume that \eqref{PtCV} is satisfied.
	Then,  for $\rho:= \lambda + \gamma/2$,
	 there exists $\ell\ge 1$
	such that:
	\begin{equation*}
		\Tsup{x\in \cX} \;\E_{x} \lp 
		\exp\lc\rho\, (\ext\wedge \tau_{\cD_\ell}) \rc \rp < \infty.
	\end{equation*}
\end{prop}

Since $(A3_F)$ is the most intricate property, 
 one may more likely suspect that it leads to restrictions.
The following proposition 
shows to what extend we may generally reply that this is not the case,
in that the convergence result actually implies a general form of property $(A3_F)$.
Note that $\mu$ in the proposition is meant to represent
 the reference measure $\zeta$ in the original statement of $(A3_F)$.
\begin{prop}
	\label{D_recA3}
	Assume that \eqref{PtCV} is satisfied.
	Then, for $\rho:= \lambda + \gamma/2$, 
	there exist $c_F >0$
	such that for any $\mu \in \Mone$,
	and any $\tZa>0$ sufficiently large,
	there exists for any $x\in \cX$
	a stopping time $U_{H}$
	satisfying the three following properties:
	\begin{align*}
		&	\PR_x(\tZa < \ext \wedge U_{H}) 
		\le \frac{c_F}{\LAg \mu \bv \heig\RAg} e^{-\rho \tZa},\qquad
		\text{ where } 
		\Lbr\ext \wedge \tZa \le \Uza \Rbr
		= \Lbr\Uza = \infty\Rbr;
		\\& \text{ and }\quad
		\E_x(X_{U_{H}} \in dy\pv U_{H} <\ext)
		\le \frac{\Ninf{\heig}}{\LAg \mu \bv \heig\RAg} 
		\E_\mu(X_{\tZa} \in dy\pv \tZa<\ext).
	\end{align*}
	Furthermore, 
	the definition of $U_{H}$
	can be extended to derive a stopping time $U_{H}^\infty$ 
	as stated in $(A3_F)$.
\end{prop}
\noindent
To deduce that $(A3_F)$ holds for any measure $\mu$ 
satisfying $\LAg \mu\bv h\RAg>0$,
one has simply to remark 
that the following inequality holds 
with $\hat{\rho} := \lambda + \gamma/4$
for any $\eps>0$
and any $\tZa$ sufficiently large: 
\begin{equation*}
(c_F/\LAg \mu \bv \heig\RAg)\cdot e^{-\rho \tZa}
\le \eps\cdot e^{-\hat{\rho} \tZa},
\end{equation*}
then to take $V:=\tZa$.

The proof of this Proposition follows in Subsection \ref{D_sec_recA3} 
the one of Proposition \ref{D_recA2}.
This proof provides also
a relevant intuition 
on the role of the parameters involved in the expression of $(A3_F)$.
\begin{rem}
The derivation of the mixing property $(A1)$
from a convergence property in total variation such as \ref{PtCV}
appears much more difficult to obtain.
It indeed requires a lower-bound with respect to some measure $\zeta$
that is uniform over the initial condition in some $\cD_\ell$,
while the total variation informs about the discrepancy between two measures.
Given our additional restriction involving $T_{\cD_L}$,
a property of uniform convergence would also at least be required,
as stated in the next Theorem~\ref{D_Approx}.
\end{rem}

\paragraph{Uniformity in the localization procedure}
\label{D_sec_local}
The constants involved in the convergences are explicitly related to the ones 
in Assumptions $\mathbf{(A_F)}$ or $\mathbf{(A)}$. 
Although the specific relation is very intricate,
it implies the following corollary of approximation: 
\begin{theo}
	\label{D_Approx}
	Assume that $\mathbf{(A_F)}$ or $\mathbf{(A)}$ holds.
Then	all the preceding results 
	hold true with the same constants involved in the convergences
	when $\ext$  	is replaced 
	by $\ext^L:= \ext \wedge T_{\cD_L}$ 
	for any $L\ge 1$ sufficiently large.
	
	Let $\alpha^L,\lambda_L, h^L$ be the corresponding QSD, extinction rates and survival capacities. 
	Then as $L$ goes to infinity, 
	$\lambda_L$ converges to $\lambda$  and  $\alpha^L,h^L$ converge to $\alpha,h$  in total variation and pointwise respectively. 
	Also, we deduce $\rho\iSv = \lambda$.
\end{theo}

\section{Proofs of the general results}
\label{D_sec_Pf}
\setcounter{eq}{0}


\subsection{Proof of Theorem~\ref{D_Th:AF}}
\label{D_sec_thAf}

Recall the expression of $c_S$ in Equation \eqref{D_eqSb}.
Thanks to Assumption $(A2)$ that controls the escape time from the transitory domain,
we may define the constant $e_\T$ as follows:
\begin{equation}\label{def_eT}
	e_\T:= \Tsup{x\in \cX} \;\E_{x} \lp 
	\exp\lc\rho\, (\ext\wedge \UDc) \rc \rp.
\end{equation}
In the following,
every time we will apply Assumption \textup{$(A3_F)$}, 
we will exploit the following parameter $
	\fl:= (2\, \cp\iSB\, e_\T )^{-1}.$
For $\tZa$ the associated deterministic 
upper-bound on $\Uza$,
we define $\tp\iDC:= \tZa + \tp\iSB$.
\\

The first step consists in the following lemma:
\begin{lem}
	\label{D_lem:K0}
	Assume that inequality \eqref{D_eqSb}\, and assumption \textup{$(A3_F)$} hold
	with the above parameters. 
	Then, there exists $C_0 > 0$ such that
	the following upper-bound holds for any $x\in E$
	and any $\tp \ge \tp\iDC$:
	\begin{equation*}
		\PR_{x} (U_{H}< \tp < \ext) 
		\le C_0 \; \PR_{\alc} (\tp < \ext).
	\end{equation*}
\end{lem}

For $\fl$ sufficiently small, 
$\{U_{H}< \tp < \ext\}$ 
is meant to be the leading part of $\PR_{x} ( \tp < \ext)$. 
We will prove indeed 
that the extension of survival during the failed coupling procedure,
for a time-length $\tZa$,
then outside $E$ (before $\tau\iET^1$) is not sufficient to compensate 
the cost of the complementary event of failure, 
i.e. $\Lbr\ext \wedge \tZa \le \Uza \Rbr
= \Lbr\Uza = \infty\Rbr$.

Our idea is to distinguish the events according to the number of failures, 
and treat them inductively.
So if a first failure is observed (implying that $\UDc^1 < U_{H}^\infty \wedge \ext \wedge t$),
we start a new turn by replacing $x$ by $X(\UDc^1)$ and $t$ by $\tp - \UDc^1$. 
To combine it efficiently with Lemma~\ref{D_lem:K0}, 
the induction is in fact not stated with a dependency on the value of  $\Nstep$,
but rather on the value of the following random variable: 
\begin{equation}
	\Nstep(t):= \sup\Lbr j \ge 0\pv \UDc^j < (\tp - \tp\iDC)\wedge \ext\wedge U_{H}^\infty\Rbr,
	\label{D_Kt}
\end{equation}
where we recall from Equation~\eqref{D_tdc}
the inductive definition the sequence $(\UDc^j)$:
\begin{equation*}
	\tau\iET^{i+1}:= \inf\{s\ge \tau\iET^i +\tZa \pv X_s \in E\}
	\wedge \ext,
	\AND
	\tau\iET^0 = 0.
\end{equation*}

Thanks to the following lemma,
whose proof marks the second step,
we will finish the initialization 
with an upper-bound on the probability of the event $\{\Nstep(t) =0\}$:
\begin{lem}
\label{D_lem:tU}
Assume that assumption $(A2)$, 
inequalities \eqref{D_eqSb}\, and assumption \textup{$(A3_F)$} hold 
with $\rho> \Rsv$. 
Then, there exists $C_F>0$ such that
the following upper-bound holds for any $x\in E$ and $\tp \ge \tp\iDC$:
\begin{equation*}
	\PR_{x} (\Uza =\infty\mVg \tp \le \UDc^1\mVg \tp < \ext) 
	\le C_F \; \PR_{\alc} (\tp < \ext).
\end{equation*}
\end{lem}

The induction property will then be propagated
thanks to the following lemma,
whose proof marks the third step.
\begin{lem}
	\label{D_indK}
	Assume that assumption $(A2)$, 
	inequalities \eqref{D_eqSb}\, and assumption \textup{$(A3_F)$} hold 
	with $\rho> \Rsv$.
	If there exists $j \ge 0$ and $C_j >0$ such that
	the following upper-bound holds for any $x\in E$
	and any $\tp \ge \tp\iDC$:
\begin{equation}\label{D_Ck}
		\PR_{x} (\tp < \ext \mVg \Nstep(t) = j) \le C_j \; \PR_{\alc} (\tp < \ext),
\end{equation}	
then the inequality holds also in the next step $j+1$
as follows 
 for any $x\in E$
and any $\tp \ge \tp\iDC$:
	\begin{equation}
		\PR_{x} (\tp < \ext \mVg \Nstep(t) = j+1) 
		\le  \tfrac{C_j}{2} \; \PR_{\alc} (\tp < \ext)
	\end{equation}
\end{lem}

With these three lemmas, we will then finish the proof of Theorem \ref{D_Th:AF}.

\subsubsection{Step 1: Proof of Lemma \ref{D_lem:K0}}
The Markov property is exploited as follows in relation to \textup{$(A3_F)$}:
\begin{align*}
	&\PR_{x} (\Uza < \tp < \ext)
	= \E_{x} \left[
	\PR_{X\lc \Uza \rc } (\tp - \Uza< \wtd{\ext})\pv
	\Uza < \ext\right]\\
	&\hspace{1cm}
	\le \Cp \; \E_{\alc} \left[
	\PR_{X \lc \UCa \rc } (\tp - \tZa< \wtd{\ext})\pv
	\UCa < \ext \right]\\
	&\hspace{1cm}
	\le   \Cp \; \PR_{\alc} \left[ \tp - \tZa  < \ext \right],
\end{align*}
where we exploited  assumption \textup{$(A3_F)$} 
for the first inequality
with the fact that $\Uza\le \tZa$
holds on the event $\{\Uza < \ext\}$,
then noted that $\tp - \tZa +V \ge \tp - \tZa$
for the second inequality.
Thanks to inequality \eqref{D_eqSb}, with $\tp \ge \tZa + \tp\iSB$,
the following upper-bound is derived:
\begin{equation*}
\PR_{x} (\Uza < \tp < \ext)
	\le \Cp\cdot\cp\iSB\cdot  e^{\rho\, \tZa}\;  
	\PR_{\alc} \left[ \tp < \ext \right].
\end{equation*}
Lemma \ref{D_lem:K0} is thus satisfied 
with: 
$ C_0:= \Cp\cdot\cp\iSB \cdot e^{\rho\, \tZa}.$
\epf

\subsubsection{Step 2: Proof of Lemma \ref{D_lem:tU}}
\label{D_sec_tU}
Thanks to Equation~\eqref{def_eT},
the following upper-bound holds 
a.s. on $\Lbr \Uza = \infty,
\tZa < \ext\Rbr$:
$$\E_{X_{\tZa} } 
\lc \exp\lp \rho [\wTdc \wedge \wtd{\ext}] \rp \rc
\le e_\T.$$
Thanks to the Markov inequality, this implies the next upper-bound
a.s. on the same event:
$$\PR_{X_{\tZa} } 
(\tp - \tp\iDC - \tZa \le \wTdc
\pv \tp - \tZa < \wtd{\ext})
\le e_\T\cdot e^{ - \rho [\tp - \tp\iDC - \tZa]}.$$
In combination with the Markov property, we deduce
the next upper-bound:
\begin{equation*}
\PR_{x} (\Uza =\infty \mVg \tp - \tau\iET \le \UDc^1\mVg \tp < \ext)
	\le  e_\T \; 
	e^{ - \rho [\tp - \tp\iDC - \tZa]} 
	\PR_{x} \left[ \tZa < \ext\mVg \Uza =\infty \right].
\end{equation*}
Thanks to Inequality~\eqref{D_eqSb},
we can relate the decay $e^{-\rho t}$
to $\PR_\alc(t<\ext)$. 
Thanks to \textup{$(A3_F)$}
and recalling the definition of $\fl$ as $(2\, \cp\iSB\, e_\T )^{-1}$,
we then conclude the proof of 
Lemma \ref{D_lem:tU} where:
\begin{equation*}
	C_F:=
	\dfrac{ e_\T \;\cp\iSB\;  \fl }
	{e^{\rho\, \tp\iSB}\;   \PR_{\alc}(\tp\iSB< \ext)}
	= (2 e^{\rho\, \tp\iSB}\;   \PR _{\alc}(\tp\iSB< \ext))^{-1}>0.
	\SQ
\end{equation*}

\subsubsection{Step 3: Proof of Lemma \ref{D_indK}}
Thanks to the Markov property assumed on $\Uza^\infty$
(see assumption \textup{$(A3_F)$}) and to definition \eqref{D_Kt}:
\begin{align*}
	&\PR_{x} (\tp < \ext\mVg \Nstep(t) = j+1)
	\\ &\hspace{1cm} 
	= \E_{x} \left[
	\PR_{X\lc \UDc^1 \rc } (\tp - \UDc^1< \wtd{\ext}\mVg
	\wtd{\Nstep}(\tp - \UDc^1) = j)\pv
	\Nstep(t) \ge 1
	\right]
	\\ &\hspace{1cm}
	\le C_j\; \E_{x} \left[
	\PR_{\alc } (\tp - \UDc^1< \wtd{\ext})\pv
	\Nstep(t) \ge 1
	\right],\; 
\end{align*}
where we exploited  \eqref{D_Ck}\, and the fact that $\tp - \UDc^1 \ge \tp\iDC$ 
on the event $\{\Nstep(t) \ge 1\}$.
Then, thanks to Inequality \eqref{D_eqSb}:
\begin{equation}
	\PR_{x} (\tp < \ext\mVg \Nstep(t) = j+1)
	\le C_j\cdot \cp\iSB\cdot \PR_{\alc } (\tp < \ext) 
	\cdot 	\E_{x} \left[
	\exp(\rho \, \UDc^1);\;
	\Nstep(t) \ge 1 \right].
	\label{in_prod}
\end{equation}
We then decompose $\UDc^1$ into the sum of $\tZa$ and $\wTdc$
to arrive at the next upper-bound
on the expectation in the right-hand side:
\begin{multline}\label{eq_dec_JPU}
	\E_{x} \left[
	\exp(\rho \, \UDc^1);\;
	\Nstep(t) \ge 1 \right]\; 
\\ \le	\E_{x} \Big[
	\E_{X_{\tZa}} \Big(
	\exp(\rho \, \wTdc )\pv
	\wTdc < (\tp - \tZa- \tp\iDC) \wedge \wtd{\ext} \Big)
	\cdot \exp(\rho \tZa)
	\pv    \Uza =\infty\mVg \tZa < \ext  \Big].
\end{multline}
Now, thanks to assumption \textup{$(A2)$},
the following upper-bound holds $a.s.$ on $\Lbr  \Uza =\infty\mVg \tZa < \ext \Rbr$:
\begin{equation}
	\label{in_eT}
	\E_{X_{\tZa}} \lp
	\exp(\rho \, \wTdc)\pv
	\wTdc < (\tp - \tZa- \tp\iDC) \wedge \wtd{\ext} \rp
	\le e_\T.
\end{equation}
The next upper-bound is derived for any $x\in E$
thanks to assumption \textup{$(A3_F)$}:
\begin{equation}
	\label{in_fl}
	\E_{x} \lc \exp(\rho \tZa) \pv \Uza =\infty\rc
	\le \fl = \dfrac{1}{2\, \cp\iSB\, e_\T}.
\end{equation}
Combining the four inequalities \eqref{in_prod}, 
\eqref{eq_dec_JPU}, \eqref{in_eT} and \eqref{in_fl}
yields the following upper-bound
for any $x \in E$ and any $\tp \ge \tp\iDC$:
\begin{equation*}
	\PR_{x} ( \tp < \ext \mVg \Nstep(t) = j+1)
	\le \tfrac{C_j}{2}\; \PR_{\alc } (\tp < \ext),
\end{equation*}
which concludes the proof of Lemma \ref{D_indK}. \epf

\subsubsection{The end of the proof of Theorem \ref{D_Th:AF}}

With $\Cp:= 2\, (C_0 + C_F)$, thanks to Lemmas \ref{D_lem:K0}-3, 
we deduce that the following upper-bound holds 
for any $x \in E$, any $\tp \ge \tp\iDC$ and any $j\in \bZ_+$:
$$\PR_{x} (\tp < \ext \mVg \Nstep(t) = j) 
\le  2^{-j-1} \; \Cp \cdot \PR_{\alc} (\tp < \ext).$$
With this decomposition, we simply conclude the proof of Theorem \ref{D_Th:AF}
as follows:
\begin{equation}
\PR_{x} (\tp < \ext) 
= \Tsum{j\ge 0} \PR_{x} (\tp < \ext \mVg \Nstep(t) = j) 
\le  \Cp \cdot \PR_{\alc} (\tp < \ext).
\SQ
\end{equation}

\subsection{A more refined convergence result}
\label{D_sec_PrCV}
By Corollary 5.2.1 in \cite{AV_QSD}
and Theorem \ref{D_Th:AF},
Assumption $\mathbf{(A)}$ of \cite{AV_QSD}
is nearly implied by our Assumption $\mathbf{(A_F)}$,
except that $\medcup_{\ell\ge 1} \cD_\ell= \cX$ is no longer assumed.
We let the reader check that the proofs
given in Section 5 of  \cite{AV_QSD}  apply mutatis mutandis,
except for the following facts:
\begin{itemize}
	 \renewcommand{\labelitemi}{$\star$}
	\item As presented in several examples in \cite{AV_GS}, 
	$\alpha$ might not be the unique QSD for $\cX$.
	By adapting the argument 
	given in the following proof of Proposition \ref{p:H0}
	we can nonetheless deduce that any other QSD $\nu$ must satisfy that:
	$\nu(\medcup_\ell \cD_\ell) = 0.$
	
	\item The lower-bound on $h(x)$ are only obtained for $x\in \cD_\ell$,
	so that $h(x)$ might be equal to 0 for $x\in \cX\setminus \medcup_\ell \cD_\ell$.
	\item The reasoning on the Q-process can only be applied 
	for initial conditions on \\
	$\mathcal{H}:= \{x\in \cX\pv h(x)>0\}\supset \medcup_\ell \cD_\ell$.
	\item The uniqueness of $\beta$ as a stationary distribution for the Q-process
	holds among all distributions with support on $\mathcal{H}$
	and not necessarily on $\M_1(\cX)$.
\end{itemize}

By adapting Theorems 2.1-3 in \cite{AV_QSD},
we deduce that there exists $\heig$, $\lambda$
$\gamma, C_h >0$
such that the following inequality holds 
for any $\xi\in (0,1]$ and any $\mu \in \Mone$ satisfying $\mu(\cD_\ell) \ge \xi$:
\begin{equation}
	\frl{\tp > 0}
	\Ninf{e^{\lambda t}\PR_.(t<\ext) - \heig}
	\le C_h \; e^{-\gamma \; \tp }.
	\label{D_ECvEtaloc}
\end{equation}
Moreover, there exists also $\alpha\in \M_1(\cX)$
and a family of constants $C_\alpha(\ell,\, \xi)>0$,
defined for any $\ell \ge 1$ and $\xi \in (0,1]$, 
 such that
the following property holds true for any $\mu \in \Mone$ such that $\mu(\cD_\ell) \ge \xi$ and any $\tp > 0$:
\begin{equation}
	\|\, \mu A_{\tp}(dx)  - \alpha(dx) \, \|_{TV}
	\le C_\alpha(\ell, \xi) \; e^{-\gamma \; \tp}.
	\label{D_ECvAlloc}
\end{equation}	

Moreover, $\LAg \alpha \bv \heig\RAg = 1$.
The first condition on the $Q$-process is also directly deduced. 
It only remains to prove the convergence results as they are stated.

\subsubsection{Proof of Theorem \ref{D_AllPho}}
\label{sec_AllPho}
We recall that $\mu$ is chosen in Definition 1 such that
$\NTV{\mu} \le 1$.
Denote by $\mu_+$ (resp. $\mu_-$)
the positive (resp. negative) component of $\mu$
so that: $\mu = \mu_+-\mu_-$
and $\NTV{\mu} = 1 = \mu_+(\cX) + \mu_-(\cX)$.
Let $y \in \cD_1$ and
 define the positive measures $\hat{\mu}_+$ and $\hat{\mu}_-$ as follows:
\begin{equation*}
	\hat{\mu}_+(dx) 
	:= \frac{1}{1+\mu_+(\cX)}\, [\delta_y + \mu_+(dx)]\ge 0,
	\qquad
	\hat{\mu}_-(dx) 
	:= \frac{1}{1+\mu_-(\cX)}\, [\delta_y + \mu_-(dx)]\ge 0.
\end{equation*}
Note that $	\hat{\mu}_+(\cX) = \hat{\mu}_-(\cX) = 1,$ so that 
both $\mu_+$ and $\mu_-$ are probability measures.
In particular, this implies the following identification:
\begin{equation}\label{eq_PAmu}
	\mu_+ . (e^{-t\, \laZ}\, P_t)(dx) = 
	\langle \hat{\mu}_+ \bv \heig_t \rangle
	\cdot \hat{\mu}_+ . A_{\tp}(dx)
\end{equation}
and similarly for $\mu_-$.
These measures are also constructed so as to satisfy the two following properties:
\begin{align*}
	& \mu = [1+\mu_+(\cX)] \cdot \hat{\mu}_+ 
	-  [1+\mu_-(\cX)] \cdot \hat{\mu}_-,
	\\&
	\hat{\mu}_+(\cD_1)\wedge \hat{\mu}_-(\cD_1)
	\ge  \frac{1}{1+\NTV{\mu}}\ge  1/2.
\end{align*}
Thus, first thanks to \eqref{D_ECvAlloc}, then to \eqref{D_ECvEtaloc},
there exists $\gamma$, $C>0$ independent from $\mu$
such that
the four following upper-bounds hold 
on $\mu_+$ and $\mu_-$
for any $t\ge 0$:
\begin{align*}
	\|\, \hat{\mu}_+ A_{\tp}(dx)  - \alpha(dx) \, \|_{TV}
&\le C \; e^{-\gamma \; t},
	\quad 
|\, \hat{\mu}_- A_{\tp}(dx) -\alpha(dx) \, \|_{TV}
	\le C \; e^{-\gamma \; t},
	\\|\langle \hat{\mu}_+ \bv \heig_t -\heig \rangle |
 &\le C \; e^{-\gamma \; t},
\hcm{2}  |\langle \hat{\mu}_- \bv \heig_t -\heig \rangle |
	\le C \; e^{-\gamma \; t}.
\end{align*}
The linear decomposition of $\mu$ between $\mu_+$ and $\mu_-$ 
thus implies the next estimations,
thanks also to Equation~\eqref{eq_PAmu} and to the fact that $\heig$ and the family $\heig_t$
are uniformly bounded:
\begin{align*}
	\mu . (e^{-t\, \laZ}\, P_t)(dx)
	&= [1+\mu_+(\cX)] \cdot\langle \hat{\mu}_+ \bv \heig_t \rangle
	\cdot \hat{\mu}_+ . A_{\tp}(dx)
	- [1+\mu_-(\cX)] \cdot \langle \hat{\mu}_- \bv \heig_t \rangle
	\cdot \hat{\mu}_- . A_{\tp}(dx)
	\\&
	= \Big( [1+\mu_+(\cX)] \cdot \langle \hat{\mu}_+ \bv \heig \rangle
	- [1+\mu_-(\cX)] \cdot \langle \hat{\mu}_- \bv \heig \rangle\Big)
	\cdot \alpha(dx)
	+ O_{TV}( e^{-t\, \gamma})
	\\&
	= \langle \mu \bv \heig \rangle\, \alpha(dx)
	+ O_{TV}( e^{-t\,\gamma}).
\end{align*}
Thus, we conclude the proof of Theorem \ref{D_AllPho} 
in that there exists some $C'>0$
such that the following upper-bound holds
for any signed measure such that 
$\NTV{\mu} \le 1$, 
and any $t\ge 0$:
\begin{equation*}
\NTV{\mu .  (e^{-t\, \laZ}\, P_t) 
-  \langle \mu \bv \heig \rangle\, \alpha} 
\le C'\, \exp[-t\, \gamma].\SQ
\end{equation*}

\subsubsection{Proof of Corollary \ref{D_CVAl}}

Let $\bar{\mu} = \mu - \alpha$. We recall that by definition, 
$\LAg \alpha\bv \heig_t\RAg = \LAg \alpha\bv \heig\RAg = 1$.
Then we can express the difference between $\mu A_t$ and $\alpha$
as follows in terms of $\bar{\mu}$
(recall also the analogous of Equation~\eqref{eq_PAmu} for $\mu$)
\begin{align*}
	\mu A_t - \alpha 
	&= \dfrac{\exp[t \laZ]\, \mu P_t - \LAg \mu\bv \heig\RAg\, \alpha}
	{\LAg \mu\bv \heig\RAg}
	+ \dfrac{\LAg \mu\bv \heig - \heig_t\RAg}
	{\LAg \mu\bv \heig\RAg \cdot \LAg \mu\bv \heig_t\RAg}\,
	\exp[t \laZ]\, \mu P_t
	\\&	=  \dfrac{\exp[t \laZ]\, \bar{\mu} P_t 
		- \LAg \bar{\mu}\bv \heig\RAg\, \alpha}
	{\LAg \mu\bv \heig\RAg}
	+ \dfrac{\LAg \bar{\mu}\bv \heig - \heig_t\RAg}{\LAg \mu\bv \heig\RAg }
	\, \mu A_t.
\end{align*}
Thanks to Theorem \ref{D_AllPho}, this immediately implies 
the estimate on the convergence to $\alpha$.
\epf

\subsubsection{Proof of Proposition \ref{p:H0}}

From $(H_0)$, we consider $H\subset \cX$, $t, c>0$ and $\ell\ge 1$ such that
the following lower-bound holds for any $x\in H$:
\begin{equation}
\PR_x(\tau_{\cD_\ell} < t \wedge \ext) \ge c.
\label{ineq_lb_H}
\end{equation}
Recalling that the  proof in \cite{AV_QSD}
ensures that $h$ is bounded away from zero by a positive constant
on any $\cD_\ell$, 
let  $h_\ell$ be such lower-bound.
The property of $h$ being an eigenfunction 
of the semi-group $(P_t)$
can be rephrased by saying that
$(h(X_t)e^{\rho_0 t}\idc{t<\ext})$
is a martingale.
We also recall that $h$ is non-negative.
The following upper-bound for any $x\in H$
is derived thanks to the martingale property, then the definition of $h_\ell$
and finally Inequality~\eqref{ineq_lb_H}:
\begin{align*}
h(x) 
&= \E_x\lc h(X(\tau_{\cD_\ell} \wedge t)) 
\exp[\rho_0 (\tau_{\cD_\ell} \wedge t)]\pv 
\tau_{\cD_\ell} \wedge t < \ext\rc
\\&
\ge h_\ell \PR_x(\tau_{\cD_\ell} < t\wedge \ext)
\ge c\cdot h_\ell >0.
\end{align*}
This proves the uniform lower-bound of $h$ on $H$.

For the second point,
recalling that the sequence of stopping times $\tau_{\cD_\ell}$
is necessarily decreasing, 
we deduce the following inclusion, 
that leads to the first desired inclusion:
\begin{align*}
 \Lbr x \in \cX\pv \Ex{\ell\ge 1}\PR_x(\tau_{\cD_\ell} < \ext)>0\Rbr
 &= \medcup_{\ell \ge 1}
  \Lbr x \in \cX\pv \PR_x(\tau_{\cD_\ell} <\ell\wedge  \ext)\ge 1/\ell\Rbr
  \\  &\subset \{x\in \cX; h(x)>0\}.
\end{align*}
For the reciprocal inclusion, let any $x$ be such that $h(x)>0$.
Thanks to Corollary \ref{D_CVAl}, $\delta_x A_t$ converges to $\alpha$.
Choosing $\ell\ge 1$ 
by Theorem  \ref{D_AllPho} such that $\alpha(\cD_\ell)>0$,
it implies that $\delta_x A_t(\cD_\ell)>0$ 
for $t$ sufficiently large,
thus $\PR_x(\tau_{\cD_\ell} < (t+1)\wedge \ext)>0$.
 This ends the proof of Proposition \ref{p:H0}.
\epf

\subsubsection{Proof of $(ii)$ in Corollary \ref{D_QECV}}

Assume that $\mu \in \Mone$ 
satisfies $\LAg \mu\bv 1/\heig\RAg < \infty$.
We may define $\nu(dx) 
:= \dfrac{\mu(dx)}{\heig(x)\, 
	\langle\mu\bv 1/\heig\rangle}$,
which trivially satisfies that 
$\nu \in \Mone$ and that
$\nu\cdot B[\heig]= \mu$.
%
Let $\bar{\nu} = \nu - \alpha$.
The difference between $\nu\cdot B[\heig]\cdot Q_t$ and $\beta$ is expressed as follows
in terms of $\bar \nu$ and of a real-valued bounded measurable function $f$:
\begin{align*}
	\LAg \nu\cdot B[\heig]\cdot Q_t\bv f\RAg - \LAg\beta\bv f\RAg
	&= \dfrac{ \LAg \nu\bv e^{t \laZ} P_t\bv \heig\cdot f\RAg 
		- \LAg \nu\bv \heig\RAg\cdot \LAg \alpha\bv \heig\cdot f\RAg}
	{\LAg \nu\bv \heig\RAg}
	\\&
	=   \dfrac{ \LAg \bar{\nu}\bv e^{t \laZ} P_t\bv \heig\cdot f\RAg 
		- \LAg \bar{\nu}\bv \heig\RAg\cdot \LAg \alpha\bv \heig\cdot f\RAg}
	{\LAg \nu\bv \heig\RAg}
	\\&
	\le \dfrac{\|\bar{\nu}\|_{TV}}
	{\LAg \nu\bv \heig \RAg} \cdot \Ninf{\heig} \cdot \Ninf{f}\cdot e^{-\gamma t},
\end{align*}
thanks to Theorem \ref{D_AllPho}.
Moreover, by the definition of $\beta$ and $\nu$:
\begin{align*}
	&\dfrac{\|\nu - \alpha\|_{TV}}
	{\LAg \nu\bv \heig \RAg} 
	= \left\|\dfrac{\mu}{\langle\mu\bv 1/\heig\rangle} 
	- \beta \right\|_{1/\heig}
	\cdot \dfrac{\langle\mu\bv 1/\heig\rangle}
	{\LAg \mu\bv \idg{} \RAg} 
	= \|\mu - \langle\mu\bv 1/\heig\rangle\, \beta\|_{1/\heig}.
\end{align*}
Injecting this equality into the preceding upper-bound
yields $(ii)$ and conclude the proof of Corollary \ref{D_QECV}.
\epf

\subsection{Proof of Propositions \ref{D_recA2}}
For $x\in \cX$ and $t\ge 0$, 
define $\nu_{x, t}(dy):= (e^{\lambda t} \delta_x P_t -\heig(x) \alpha)_+(dy) \ge 0$.
Thanks to \eqref{PtCV}, we know that there exists $C, \gamma >0$
such that $\nu_{x, t}(\cX) \le C e^{-\gamma t}$.
Let $\rho:= \lambda + \gamma/2$
and $t_E:= 2\cdot[\gamma \cdot\log(4 C)]^{-1}$.
We thus ensure that $C\cdot e^{-(\lambda + \gamma) t_E} \le (1/4)\cdot e^{-\rho t_E}.$

Since $\cX = \medcup_{\ell\ge 1} \cD_\ell$, 
we choose $\ell$ sufficiently large to ensure 
$\alpha(\cD_\ell^c) \le \frac{e^{-\gamma\, t_E/2}}{4\Ninf{\heig}}$.
Recalling the definition of $\nu_{x, t_E}$, 
and denoting as $\cD_\ell^c$ the complementary of $\cD_\ell$,
the above results imply the following upper-bound
for any $x\in \cX$:
\begin{align}
	\delta_x P_{t_E}(\cD_\ell^c) 
&\le e^{-\lambda t_E} \heig(x) \alpha(\cD_\ell^c)
	+ \nu_{x, t_E}(\cX)
\notag\\	&\le e^{-\rho t_E}/2.
	\label{D_RetC}
\end{align}
We split $\bR_+$ into time-intervals of length $t_E$,
so as to decompose the following expectation for any $K\ge 2$:
\begin{multline}
	\E_{x} \lp \exp\lc\rho\, (\ext\wedge \tau_{\cD_\ell}\wedge (K t_E)) \rc \rp
	\\
	\le e^{\rho t_E} + \sum_{k = 1}^{K-1} \exp([k+1] t_E)\, 
	\PR_x(\frl{k'\in \II{1, k}} X_{k' t_E} \notin \cD_\ell\pv  k t_E < \ext).
	\label{recA21}
\end{multline}
Thanks to the Markov property and to Inequality~\eqref{D_RetC}, 
an elementary induction proves the following upper-bound 
for any $k\ge 1$:
\begin{equation}
	\PR_x(\frl{k'\in \II{1, k}} X_{k' t_E} \notin \cD_\ell\pv  k t_E < \ext)
	\le \frac{e^{-\rho k t_E}}{2^k}.
	\label{recA22}
\end{equation}
Thanks to Inequalities~\eqref{recA22} and \eqref{recA21}, 
we obtain the following upper-bound:
\begin{equation*}
	\textstyle
	\E_{x} \lp \exp\lc\rho\, (\ext\wedge \tau_{\cD_\ell}\wedge (K t_E)) \rc \rp
	\le e^{\rho t_E} (1+ \sum_{k\ge 1} 2^{-k}) \le 2 e^{\rho t_E} < \infty.
\end{equation*}
Letting $K$ tend to infinity concludes the proof of Propositions \ref{D_recA2}.
\epf

\subsection{Proof of Propositions \ref{D_recA3}}
\label{D_sec_recA3}
Let $x\in \cX$, $\mu\in \M_1(\cX)$ and $\tZa>0$.

\paragraph{Step 1:  Definition of $\Uza$}
\textcolor{white}{:}

We define:
\begin{equation*}
\nu^{\tZa}_{x, \mu} 
:=  \lp \delta_x P_{\tZa}
- \frac{\Ninf{\heig}}{\LAg \mu \bv \heig\RAg} \mu P_{\tZa} \rp_+.
\end{equation*}
We impose that $U_{H}$ takes values $\tZa$ or $\infty$ 
in such a way that $\nu^{\tZa}_{x, \mu}$ exactly correspond to the harvested measure:
\begin{equation}
\E_x(X_{\tZa} \in dy\pv U_{H} = \tZa<\ext) 
= \delta_x P_{\tZa}(dy) - \nu^{\tZa}_{x, \mu}(dy)\ge 0.
\label{UnuS}
\end{equation}
A  natural choice of $U_{H}$
is defined through $U$ being a uniform random variable on $(0, 1)$,
independent of the process $X$.
The choice of $\nu^{\tZa}_{x, \mu}$ was made to ensure,
by the Radon-Nikodym Theorem, that 
$\frac{\partial \nu^{\tZa}_{x, \mu}}{\partial \delta_x P_{\tZa}}: 
\cX \rightarrow [0,1]$
can be  defined as a density for the measure $\delta_x P_{\tZa}$.
So we simply impose $U_{H} = \infty$
if $U \ge \frac{\partial \nu^{\tZa}_{x, \mu}}{\partial \delta_x P_{\tZa}} (X_{\tZa})$,
and $U_{H} =\tZa$ otherwise.
With this choice, \eqref{UnuS} is satisfied
and $U_{H}$ can clearly be extended into $U_{H}^\infty$
by exploiting a new uniform random variable, independent of the previous construction,
for the coupling after each failure. 

\paragraph{Step 2:  Control of the densities}
\textcolor{white}{:}

Thanks to the definitions of $\Uza$ and $\nu^{\tZa}_{x, \mu}$,
the following inequality is straightforward:
\begin{equation*}
\E_x(X_{U_{H}} \in dy\pv U_{H}<\ext) 
\le \frac{\Ninf{\heig}}{\LAg \mu \bv \heig\RAg}
\E_{\mu}(X_{\tZa} \in dy\pv \tZa<\ext).
\end{equation*}

\paragraph{Step 3:  Control of failures}
\textcolor{white}{:}

With this definition of $\Uza$, recalling \eqref{UnuS}:
\begin{equation}
\PR_x(\Uza = \infty)
= 1 - (\delta_x P_{\tZa}(\cX) - \nu^{\tZa}_{x, \mu}(\cX))
= \nu^{\tZa}_{x, \mu}(\cX),
\label{Uinf}
\end{equation}
which is thus the quantity for which we want an upper-bound.

Consider the following positive measure:
\begin{equation*}
\hat{\nu}^{\tZa}_{x, \mu}
:= e^{-\lambda \tZa} \cdot \Big[
\Big(e^{\lambda \tZa}\delta_x P_{\tZa} - \heig(x) \alpha\Big)_+
 + \frac{\Ninf{\heig}}{\LAg \mu \bv \heig\RAg}
  \cdot\Big(\LAg \mu \bv \heig\RAg \alpha - e^{\lambda \tZa}\mu P_{\tZa}\Big)_+\Big].
\end{equation*}
Thanks to Inequality~\eqref{PtCV}  with $\rho:= \lambda + \gamma/2$,
the mass of this measure can be efficiently upper-bounded:
\begin{equation}
\hat{\nu}^{\tZa}_{x, \mu}(\cX) 
\le C e^{-(\lambda +\gamma) \tZa} \Big(1 + \frac{\Ninf{\heig}}{\LAg \mu \bv \heig\RAg}\Big)
\le \frac{c_F}{\LAg \mu \bv \heig\RAg} e^{-\rho \tZa},
\label{nuSup}
\end{equation} 
where $c_F:= 2 C\cdot \Ninf{\heig}$ is independent of $x$, $\mu$ and $\tZa$.

On the other hand, $\hat{\nu}^{\tZa}_{x, \mu}$ is such that
the following property holds:
\begin{equation*}
\delta_x P_{\tZa}(dy) 
	\le \frac{\Ninf{\heig}}{\LAg \mu \bv \heig\RAg} \mu P_{\tZa}(dy)
	+ \hat{\nu}^{\tZa}_{x, \mu}(dy),
\end{equation*}
which implies that $\nu^{\tZa}_{x, \mu} \le \hat{\nu}^{\tZa}_{x, \mu}$.
Combining it with \eqref{Uinf} and \eqref{nuSup}, the intended inequality is obtained:
\begin{equation*}
\PR_x(\Uza = \infty)
\le \frac{c_F}{\LAg \mu \bv \heig\RAg} e^{-\rho \tZa}.
\end{equation*}
Since $\tZa$  is indeed allowed to take any sufficiently large values,
this concludes the proof of Proposition \ref{D_recA3}.
\epf

\subsection{Proof of Theorem \ref{D_Approx}}

Recall that we wish to describe the approximations of the previous dynamics
when extinction happens at 
$\ext^L:= \ext\wedge T_{\cD_L}$
instead of $\ext$.

There is an explicit relation 
between all the constants 
introduced in the proofs of Theorems 
\ref{D_AllPho}-4 
(requiring also the proofs in \cite{AV_QSD}).
Moreover,
the proof actually relies 
on a single value of 
$\rho > \rho_S$ and a specific set $E$.
Note that for any $L$ such that 
$E \subset \cD_L$, we have:
$$	\Tsup{x\in \cX} \;\E_{x} \lp 
\exp\lc\rho\, (\ext^L\wedge \UDc) \rc \rp 
\le 
\Tsup{x\in \cX} \;\E_{x} \lp 
\exp\lc\rho\, (\ext\wedge \UDc) \rc \rp 
:= e_\T.
$$
Likewise,
Assumption \textup{$(A3_F)$} extends naturally 
for $\ext^L$.
 \textup{$(A1)$} is stated with extinction already occurring at the exit 
 of some set $\cD_{L(\ell)}$ prescribed by the value of $\ell$.
 Considering the proof of Theorem~\ref{D_AllPho} in Subsection~\ref{sec_AllPho},
 we see that $(3.3)$ and $(3.4)$ are only required for $\ell =1$ and $\xi =1/2$.
 Once these two values are fixed,
 the proof given in \cite[Section 5]{AV_QSD}
 treats uniformly initial conditions $\mu$
such that $\mu(\cD_1) \ge 1/2$
by exploiting \textup{$(A1)$} a finite number of times.
One can thus identify an upper-bound $L_\veebar\ge 1$
of the values $L(\ell)$ 
involved in the successive applications of \textup{$(A1)$}. 
%
So it suffices to take $L$ such that 
$\cD_{L_\veebar} \subset \cD_L$
to ensure that all the results extend
for the extinction time $\ext^L$ 
instead of $\ext$.
Under this condition, 
our proof ensures that 
uniform exponential quasi-stationary convergence
also holds for the process 
with extinction at time $\ext^L$ 
and that the constants involved in the convergences can be taken uniformly 
over these values $L$. 
\\

To compare $\lambda$ to $\lambda_L$,
we can observe that for any $t>0$:
$$\frac{-1}{t} \log \PR_\alc(t<\ext)
\le \frac{-1}{t} \log \PR_\alc(t<\ext^L)$$
so that $\lambda \le \lambda_L$
is deduced by taking the limit $t\ifty$ 
and exploiting the convergence to the survival capacities.
The same argument ensures that $\laZ_L$ is a decreasing sequence in $L$.

Assume then by contradiction that there exists 
$\eta>0$ such that 
$\lim_L \laZ_L \ge \laZ + \eta$.
Recall 
that we have an explicit upper-bound 
$\|\heig_*\|_\infty$ that is valid uniformly for the functions $\heig^L$ and $\heig$.
Thanks to Theorem \ref{D_AllPho}
and the analogous result with $\ext^L$,
for $t$ sufficiently large,
the next two  estimates on the survival probability holds:
\begin{equation*}
	e^{\laZ t} \PR_\alc(t< \ext)
	\ge \frac{1}{2}\LAg \alc\bv \heig\RAg,
	\qquad
	e^{\laZ_L t} \PR_\alc(t< \ext^L)
	\le 2 \|\heig_*\|_\infty.
	\label{D_lbsv}
\end{equation*}
By exploiting the property of $\eta$, 
we deduce:
\begin{align}
	0<\LAg \alc\bv \heig\RAg
	&\le 2 e^{(\laZ - \laZ_L) t} \cdot e^{\laZ_L t}\cdot \PR_\alc(t< \ext^L)
	+ 2 e^{\lambda t}\ 
	|\PR_\alc(t< \ext) - \PR_\alc(t< \ext^L)|
\notag	\\ &\le 4\, \|\heig_*\|_\infty e^{-\eta t}
	+ 2 e^{\lambda t}\ 
	|\PR_\alc(t< \ext) - \PR_\alc(t< \ext^L)|.
	\label{D_ubsv}
\end{align}
The first term in the upper-bound becomes negligible uniformly over $L$ by taking $t$ sufficiently large.
In order to obtain a contradiction,
we merely have to prove that 
$\PR_\alc(t< \ext^L)$ converges to $\PR_\alc(t< \ext)$ 
at any fixed (large) time $t$. The difference is 
$\PR_\alc(T_{\cD_L}< t< \ext)$
thus upper-bounded by $\PR_\alc(T_{\cD_L}< t)$.
So our first aim is to prove that a.s. 
$\lim_L T_{\cD_L} = \infty$.

Assume by contradiction 
that the limit $T_\infty$ 
of this increasing sequence 
is at a finite value.
Yet, 
thanks to $(\underline{A0})$ and to the fact that $X$ is càd-làg,
$X_{T_\infty-}\in \cX$ belongs to $\cD_M$ for some $M$.
Thus, there exists a vicinity to the left 
of $T_\infty$ on which
$T_{\cD_L}$ for $L>M$ cannot happen.
But this precisely contradicts the definition of $T_\infty$. 
Consequently,  
$\lim_L T_{\cD_L} = \infty$
holds a.s.
and $\PR_\alc(t< \ext^L)$ converges to 
$\PR_\alc(t< \ext)$ as $L\ifty$.
The contradiction with \eqref{D_ubsv}
makes us conclude that 
$\laZ_L$ tends to $\laZ$ as $L\ifty$.
\\

The next step is to look at the survival capacities,
thanks again to Theorem \ref{D_AllPho}
(with the measures evaluated on $\cX$).
The following upper-bound holds for any $x\in \cX$:
\begin{equation*}
	|\heig(x) - \heig_L(x)|
	\le e^{\laZ t} |\PR_x(t< \ext) - \PR_x(t< \ext^L)|
	+ |e^{\laZ t} - e^{\laZ_L t}|
	+ C \, e^{-\gamma t}.
\end{equation*}
Again, we can choose $t$ sufficiently large
to make $C \, e^{-\gamma t}$ negligible.
We already know that $\laZ_L$ tends to $\laZ$
and as previously, 
we prove that $\PR_x(t< \ext^L)$
tends to $\PR_x(t< \ext)$,
as $L\ifty$.
This concludes the punctual convergence of
$\heig_L$ to $\heig$.
The conclusion would be the same 
if one replaces
$x$ by any probability measure $\mu$, 
for instance $\alpha$.

Concerning the QSD:
\begin{align*}
	&\NTV{\alpha - \alpha_L}
	\le \NTV{e^{\laZ_L t}\delta_\alpha P_t^L 
		- \LAg \alpha \bv \heig_L\RAg \alpha_L}
	+ |e^{\laZ_L t} - e^{\laZ t}|
	\\&\hcm{4}
	+ |\LAg \alpha\bv \heig_L -\heig\RAg|
	+ e^{\laZ t} \NTV{ \delta_\alpha P_t^L 
		-  \delta_\alpha P_t},
\end{align*}
where as $L \ifty,$ for $t$ fixed, 
we have just shown that the following quantity tends to 0
by proving that $\lim_L T_{\cD_L} = \infty$
holds a.s.:
\begin{equation*}
		\NTV{ \delta_\alpha P_t^L 
		-  \delta_\alpha P_t}
	= \PR_\alpha (T_{\cD_L}<t< \ext)
	\rightarrow 0.
\end{equation*}
Thanks again to Theorem \ref{D_AllPho}
and the previous convergence results
to $\alpha_L$, $\heig$ and $\lambda$,
the  right-hand side
can be made negligible by taking $t$ then $L$ sufficiently large,
concluding the convergence of $\alpha_L$
to $\alpha$ in total variation.
This concludes the proof of Theorem \ref{D_Approx}.
\epf

\section*{Applications}
\label{D_sec_Appli}
\setcounter{eq}{0}

 \section{Mutations compensating a drift leading to extinction}
 \label{D_sec_Adapt3}
 
 
 \subsection{A first simple process}
 \label{D_MComp}
 
 We recall that we wish to prove 
 uniform exponential quasi-stationary convergence 
 for the following process~:
 \begin{equation}
X_t 
 = x - v\, t\, \mathbf{e_1}
 + \Tsum{i\le N_t} W_i,
 \label{D_AdS3}
 \end{equation}
 with a state-dependent extinction rate 
given by $\rho_e: x\mapsto \|x\|^2$.
 The number $N_t$  of mutations  at time $t$
 is given as a classical Poisson process on $\bR_+$.
 Each mutation effect $W_i$ is distributed as a normal variable with covariance matrix $\sigma^2 I_d$, and drawn independently of each others and of $N_t$.
 Between jumps, the process is translated at constant speed $v>0$
 along the first coordinate (i.e. along $\mathbf{e_1}$).

\begin{theo}
	\label{D_prop.A3}
	Consider $P$ the semi-group associated to the process $X$ 
as above (including the extinction).
Then, for any $v, \sigma>0$,	
$P$ displays a uniform exponential quasi-stationary convergence
	with some characteristics
	$(\alpha, h, \lambda) \in \M_1(\bR^d)\ltm B(\bR^d)\ltm \bR_+$
	(cf Definition \ref{UEQS}).
	Moreover, $h$ is positive and bounded.
\end{theo}

\subsection{The main required properties}
\label{D_sec_Mrp}

This application 
is related to non-local reaction-diffusion equations 
with a drift term.
The one dimensional case
has been studied recently by \cite{CH18},
with existence results obtained with compactness argument,
and in Section 2 of \cite{CG20},
with the use of Lyapunov functions.

To highlight the generality of our approach,
we specify next the main properties of $X$
that we exploit.
We consider generally a càd-làg process $X$ on $\bR^d$,
confronted to an extinction 
at a state-dependent rate 
given by $\rho_e:\bR^d\mapsto \bR_+$, and of the following form:
\begin{equation}
X_t 
= x - v\, t\, \mathbf{e_1}
+ \int_{[0,t] \times \bR^d \times \bR_+ } 
w \; \idc{u\le g(X_{s^-},w )}\, M(ds,dw,du).
\label{D_eqS}
\end{equation}
where 
$M$ is a Poisson Random Measure  (PRaMe) 
over $\bR_+\times\bR^d\times\bR_+ $, 
with intensity 
$\pi (ds,dw,du) = ds \; dw \; du  $, 
while $g(x, w)$ describes the jump rate from $x$ to $x+w$.
In our focal example, 
$$g(x, w):= (\dfrac{1}{\sqrt{2\pi}\sigma})^d
\exp\Big( -\frac{\|w\|^2}{2 \sigma^2}\Big).$$

 The infinitesimal generator $\cL$ of such generic process 
is defined on all $C^1$ and bounded function $f$ on $\bR^d$ as follows:
\begin{equation*}
	\cL f(x) 
	:= -v \partial_{x_1}f(x) + \int_{\bR^d} (f(x+w) - f(x))\cdot  g(x, w) dw
	- \rho_e(x) f(x).
\end{equation*}
The dynamics prescribed by the dual $\cL^\star$ of $\cL$ by the equation $\partial_t u = \cL^\star u$,
which is the starting point of  \cite{CH18}, then corresponds to the dynamics in time of the density of the measure-valued process $\mu P_t$ (see notably Section 2 in \cite{CG20}).

The next properties are stated 
for these measurable functions $g$ and $\rho_e$,
with $B(x, r)$ the open ball around $x$ of radius $r$ for the Euclidian norm.
The fact that an upper-bound holds locally in $x$ means 
that for any compact subset $K$ of $\cX$,
the upper-bound holds uniformly fpr $x\in K$.
\\

\noindent 
\textbf{Assumption $(\mathbf{P})$}(for Piecewise-Deterministic):
\begin{itemize}
\item[$(\mathbf{P}1)$] $\rho_e$ is locally upper-bounded and 
$\TlimInf{\|x\|}\rho_e(x)= +\infty$.

Also, explosion implies extinction:
$\ext \le \Tsup{\ell\ge 1} T_{\cD_\ell}$.

\item[$(\mathbf{P}2)$] The jump-rate
$\rho_J(x):= \int_{\bR^d} g(x, w) \, dw$
is locally upper-bounded.

\item[$(\mathbf{P}3)$] Locally in $x$, there exists $0<\dS<\aS$
such that the restriction of $g$ 
to $\cX\times B(\aS\cdot \mathbf{e_1}, \dS)$
is  lower-bounded.

\item[$(\mathbf{P}4)$] The jump size is  tight locally in $x$.

\item[$(\mathbf{P}5)$] The density for each jump vector $w$
is  upper-bounded locally in $x$.
\end{itemize}

\begin{theo}
\label{D_th.A3}
Provided the above conditions  $(\mathbf{P})$ are satisfied,
	$P$ displays a uniform exponential quasi-stationary convergence
with some characteristics
$(\alpha, h, \lambda) \in \M_1(\bR^d)\ltm B(\bR^d)\ltm \bR_+$.
Moreover, $h$ is positive an bounded.

Besides, the Q-process exists 
and is exponentially ergodic with weight $1/h$
as stated in Corollary \ref{D_QECV},
while the uniformity in the localization procedure
holds as stated in Theorem~\ref{D_Approx} 
for $\cD_{\ell}:= B(0, \ell \aS)$.
\end{theo}

The proof of Theorem \ref{D_th.A3} is given in Subsection \ref{D_sec_th.A3}.
It entails the proof of Theorem~\ref{D_prop.A3}
since it is elementary that 
the process described in \eqref{D_AdS3}
satisfies $(\mathbf{P})$.
Let us nonetheless first clarify the meaning of these different properties.
\\

To fix ideas, 
let us consider any compact $K$ subset of $\cX$.
Thanks to Property $(\mathbf P1)$ 
there exists $\rho_\vee$ such that $\rho_e(x) \le \rho_\vee$
for any $x\in K$. 
This makes it possible to have simple lower-bounds on the survival of a given trajectory of $(X_t)_{t\ge 0}$
provided it remains confined in the compact $K$.
On the other hand, 
the fact that $\rho_e(x)$ tends to infnity as $x$ tends to infinity 
makes it possible to justify the complementary of $B(0, \ell)$ as transitory for $\ell$ sufficiently large
(according to property $(A2)$).

Thanks to Property $(\mathbf P2)$,
there exists $\rho_J^\vee>0$ such that 
for any $x\in K$, $\rho_J(x) \le \rho_J^\vee$.
We can thus consider events of positive probabilities
such that the corresponding trajectories of $X$
have no other jumps than the one we carefully describe.

According to Property $(\mathbf P3)$,
there exists also $g_\wedge>0$ and $0<\dS<\aS$ 
such that the following inequality holds for any $x\in K$
and any $w\in B(\aS\, \mathbf{e_1}, \dS)$, $g(x, w)\ge  g_\wedge$.
This will make it possible to consider trajectories 
in which each jump compensate the drift of the previous time-interval,
with a small variation.

In addition, according to Property $(\mathbf P4)$,
for any $\eps >0$, there exists $w_{\vee}$ such that
the following upper-bound holds uniformly in $x\in K$: 
$\textstyle \int_{\bR^d} g(x, w)\,\idc{\|w\|\ge w_{\vee}}\, dw \le  \eps\cdot \rho_J(x).$
This makes it possible to restrict the size of the jumps with a probability close to 1.

Finally, according to Property $(\mathbf P5)$,
there exists  $g_\vee$ 
such that the following upper-bound holds uniformly in $x\in K$ and $w\in \bR^d$:
$g(x,w) \le g_\vee\cdot \rho_J(x).$
Thanks to this property, we will deduce upper-bounds of the marginal density of $X_t$
according to the Lebesgue measure after some jumps.

\subsection{Proof of Theorem \ref{D_th.A3}}
\label{D_sec_th.A3}

We aim at proving Assumption $\mathbf{(A_F)}$
for the sequence $\cD_{\ell}:= B(0, \ell \cdot \aS)$.

$(\overline{A0})$ is clearly satisfied.
The proof of \textup{$(A1)$} is deduced from
the following proposition, whose  proof
is deferred to the end of this subsection:
\begin{prop}
\label{D_lem.Ad3.2}
Under $(\mathbf{P}1,2,3)$, for any $\ell\ge 1$,
with $L:= \ell+2$,
  there exists $c, t>0$ such that:
  \begin{equation*}
\frl{x \in \cD_{\ell}}
  \hspace{.5cm}
  \PR_x \lc {X}_{\tp}\in dx \pv
  \tp < \ext \wedge T_{\cD_L} \rc 
  	\ge \cp\; \idg{\cD_\ell}(dx).	
  \end{equation*}
\end{prop}
In particular, 
it implies that Assumption $(A1)$ holds with $\alc$ 
uniform over $\cD_1$.
From Lemma 3.0.2 in \cite{AV_QSD}, 
noting that $\alc(\cD_1)>0$ in particular,
we can (explicitly) deduce a strict upper-bound $\rho$ of $\rho_S$.
Since the extinction rate outside of $\cD_\ell$ 
tends to infinity while $\ell\rightarrow \infty$,
for any $\rho>0$, we can find some $L\ge 1$ 
such that  assumption $(A2)$ holds true for $E:= \cD_\ell$ 
(cf. Subsection 4.1.2 in \cite{AV_QSD}).
The proof of Assumption $(A3_F)$
for these choices
is a clear consequence of the next proposition,
whose proof
is given in the next  Subsection \ref{D_sec_AEF}: 
\begin{prop}
\label{D_prop:AEF}
Suppose that Assumption $(\mathbf{P})$ holds true.
Consider any $\rho\iET > \Rsv$
and $\ell\iET\ge 1$
such that the set $E = \bar{B}(0, \ell\iET)$ 
satisfies $\frl{y\notin E} \rho_e(y) \ge \rho\iET.$
Set also any $\eps>0$.
Then, there exists $\tZa, t_V, c>0$,
such that for any $x\in E$, there exists
a stopping time $U_{H}$
with the following properties:
\begin{align}
\PR_x (\Uza = \infty,\; \tZa < \ext) 
\le \fl \cdot \exp[\Rsv\, \tZa],\qquad
& \text{ where } \Lbr\ext \wedge \tZa \le \Uza \Rbr
= \{\Uza = \infty\},
\label{D_Ad3.N}
\\ 
\PR_x ( X(\Uza) \in dx' 
\pv \Uza < \ext)
&\le  \cp \,\PR_{\alc} \big(X(t_V) \in dx'
\pv t_V < \ext\big),\notag
\end{align}
including the fact that
its definition can be extended into the one of some $U_{H}^\infty$
as specified in Assumption $(A3_F)$.
\end{prop}

Given these two propositions,
we can conclude that Assumption $\mathbf{(A_F)}$ holds true.
By Theorem~\ref{D_AllPho}, Corollary \ref{D_QECV} and Theorem~\ref{D_Approx} 
it directly implies Theorem \ref{D_th.A3}.

\subsubsection{Proof of Proposition  \ref{D_prop:AEF}} 
\label{D_sec_AEF}

With the notations of the proposition, 
we first define $\tZa$ by the relation:
\begin{equation*}
\exp[\Rsv\cdot (2\ \ell\iET / v) -  (\rho\iET-\Rsv)\cdot(\tZa - 2\ \ell\iET / v)]
= \fl /2.
\end{equation*}
The left-hand side is decreasing and converges to 0
 when $\tZa \rightarrow \infty$, so that
$\tZa$ is well-defined.
Let $T_{J}$ be the first jump time of $X$.
On the event $\Lbr \tZa < T_{J}\Rbr$,
 we set $\Uza = \infty$.
The choice of $\tZa$ is done to ensure 
that the probability associated to the failure
 is indeed exceptional enough 
 (with threshold $\fl/2$ and time-penalty $\Rsv$).
Any jump occurring before $\tZa$ occurs from a position
$X(T_{J}-) \in \bar{B}(0, \ell\iET + v\, \tZa):= K$.
Thanks to Assumption $(\mathbf{P4})$, 
we can then define $w_{\vee}$ such that:
\begin{equation*}
\frlq{x\in K} 
\textstyle \int_{\bR^d} g(x, w) \,\idc{\|w\|\ge w_{\vee}}\, dw 
\le  \fl/2\cdot \exp[\Rsv\, \tZa].
\end{equation*}
A jump size larger than $w_{\vee}$ is then the other criterion of failure.

On the event  $\Lbr  T_{J} \le \tZa\Rbr
\cap \Lbr T_{J} < \ext\Rbr
\cap \Lbr \|W\| \le w_{\vee}\Rbr$,
where  $W$ is the size of the first jump 
(at time $T_{J}$),
we thus set $\Uza:= T_{J} \le \tZa$.
Otherwise $\Uza:= \infty$.

The proof that $\Uza$ extends to some $U_{H}^\infty$ 
as stated in $(A3_F)$
is elementary and the reader will be spared these details.
Let us just say
that $U_{H}^\infty$ is defined on the event
$\{U_{H}^\infty < \ext\}$, 
as the first jump time $T'_J$ to satisfy the three following conditions:
$T'_J$ is the first jump time of $X$ 
after $\tau\iET^i$ for some $i\ge 0$
(cf  \eqref{D_tdc});
$T'_J \le  \tau\iET^i + \tZa$ for this value of $i$;
and $\|\Delta X_{T'_J}\|\le w_\vee$.

In particular, $ \Lbr\ext \wedge \tZa \le \Uza \Rbr
= \Lbr\Uza = \infty\Rbr$ is clearly satisfied.
\\

We prove next that the failures are indeed exceptional enough:
\begin{equation*}
\PR_x (\Uza = \infty,\; \tZa < \ext)
\le \PR_x ( \tZa \le T_{J}\wedge\ext)
+ \PR_x ( T_{J}<\ext\wedge\tZa\mVg \|W\| > w_{\vee}).
\end{equation*}
By the definition of $w_{\vee}$, we deal with the second term:
\begin{align*}
\PR_x ( T_{J}<\ext\wedge\tZa\mVg \|W\| > w_{\vee})
&\le \PR_x (\|W\| > w_{\vee} \bv T_{J}<\ext\wedge\tZa) 
\\&\le  \fl/2\cdot \exp[\Rsv\, \tZa].
\end{align*}
On the event $\Lbr \tZa \le T_{J}\wedge\ext \Rbr$
it holds a.s. for any $t\le \tZa$ that 
$X_t = x - v\, t\, \mathbf{e_1}$.
Thus, $X$ is outside of $E$ in the time-interval 
$[(2\ell\iET/v), \tZa]$, with an extinction rate at least $\rho\iET$.
By the definition of $\tZa$, it implies:
$$\PR_x ( \tZa \le T_{J}\wedge\ext)
\le \exp[-  \rho\iET\ (\tZa - 2\ \ell\iET / v)]
\le \fl/2\cdot \exp[\Rsv\, \tZa].$$
This concludes \eqref{D_Ad3.N}.
\\

On the other hand, 
recall that a.s. on the event $\{U_H <\infty\}$,
$X_{U_H} = W + X(T_{J}-)$
where $X(T_{J}-) \in K$.
Thus, thanks to Assumption $(\mathbf{P5})$,
there exists $g_\vee>0$
such that the following upper-bound on the density 
of $X(\Uza)$
holds uniformly in $x\in E$:
 \begin{equation*}
 \PR_x ( X(\Uza) \in dx' 
 \pv T_{J}< \tZa\wedge \ext\mVg \|W\| \le w_{\vee})
 \le g_\vee\, \idc{x'\in \bar{B}(0, \ell\iET + v\, \tZa+w_{\vee})} dx'.
 \end{equation*}
We know also from Proposition \ref{D_lem.Ad3.2} that there exists $t\iMix, c\iMix>0$ 
such that:
\begin{equation*}
\PR_{\alc} \big(X(t\iMix) \in dx\pv t\iMix < \ext\big)
\ge c\iMix\, \idc{x'\in \bar{B}(0, \ell\iET + v\, \tZa+w_{\vee})} dx'.
\end{equation*}
With $t:= t\iMix$, $c:= g_\vee/ c\iMix$
 and thanks to Inequality \eqref{D_Ad3.N}, this concludes the proof of Proposition \ref{D_prop:AEF}.
\epf

\subsubsection{Proof of Proposition \ref{D_lem.Ad3.2}}

%
%
We consider a characteristic length 
of dispersion given by $r:= \dS/4$.
Given some $\ell\ge 1$, $x_I\in \cD_\ell$, $L:= \ell+2$ and $c>0$,
we propose the following definition of the range of triplets time/density/position
that can be reasonably reached by the process,
starting from the vicinity of $x_I\in \cD_\ell$
and restricted on specific domains (of the form $\cD_L$):
 \begin{align*}
 &\cR^{(L)}:= \Lbr (t, c, x_F)\in \bR_+\times\bR_+\times \bR^d \pv
 \frl{x_0 \in B(x_I, r)}
\right.
 \\&\hcm{2.5}
 \left. 
 \PR_{x_0} \lp X_t\in dx
 \pv t < T_{\cD_L}\wedge \ext \rp
 \ge c\; \idg{B(x_F, r)}(x)\,  dx\Rbr.
 \label{D_cR}
 \end{align*}

The proof then relies on the following three elementary lemmas.
As a first step formalized in the next Lemma~\ref{D_lem:Iloc},
we show that $\cR^{(L)}$ 
contains the product of a non-empty time-interval
and of a vicinity of $x_I$,
for a sufficiently small density factor $c_0$. 
\begin{lem}
\label{D_lem:Iloc}
Given any $L\ge 3$, 
there exists 
$c_0, t_0, \delta>0$, such that
the following inclusion holds for any $x_I \in \cD_\ell$ with $\ell = L-2$:
$$[t_0, t_0 + \delta] \times [0, c_0]\times B(x_I, r) \subset \cR^{(L)}.$$
\end{lem}
The lemma is proved in Subsection \ref{D_sec_Iloc} by compensating the drift component
with exactly one jump and adjusting the time to the allowed variations in jump vector.
Thanks to the lemma, we may start to consider trajectories from any initial position in $\cD_\ell$. 
In addition, a time-interval is considered
in order to adjust as a last step the durations of the trajectories 
that link the various reference points (corresponding to both $x_I$ and $x_F$).
As the second step formalized in the next Lemma~\ref{D_lem:Iloc},
we show that we may expand the range in the feature dimension $\cX$,
at the expense of a specific increase of time and reduction of the reference density:
\begin{lem}
\label{D_lem:Rloc}
For any $L\ge 3$, with $\ell = L-2$,
there exists $t_a, c_a>0$ such that
the following implication holds uniformly
for any $x_I, x \in \cD_\ell$ and any $t, c>0$:
\begin{equation*}
(t, c, x) \in \cR^{(L)}
\quad \imp \quad
\{t+t_a\} \times \{c\cdot c_a\}\times B(x, r)\subset \cR^{(L)}.
\end{equation*}
\end{lem}
The lemma is proved in Subsection \ref{sec_lem_Rloc}, 
once more with the compensation of the drift by one jump
with enough flexibility.
This lemma is exploited inductively in the final step, cf Subsection~\ref{sec_expand},
until the whole $\cD_\ell$
belongs to the feature component of the range (for a large enough time and a small enough reference density).
We also need to justify that 
we can make the associated durations coincide,
up to a reduction of the density factor,
which exploits the next Lemma in combination 
with Lemma \ref{D_lem:Iloc}.
\begin{lem}
	\label{D_lem:PropA}
	There exists $c_P>0$ (only depending on $r$ and $d$)
such that the following implication holds for any $L\ge 3$,
for	any $x_I \in \cD_\ell$ with $\ell = L - 2$,
any $t, c>0$ such that $(t, c, x_I) \in \cR^{(L)}$, 
any $t', c'>0$ and any $x'\in \cD_L$:
	\begin{equation*}
	(t',c',x')\in \cR^{(L)} 
	\imp (t+t', c\cdot c_P\cdot c', x') \in \cR^{(L)}.
	\end{equation*}
\end{lem}
The lemma is proved as the third step in Subsection \ref{sec_lem_Rloc}. 
The forth step concludes the proof of Proposition \ref{D_lem.Ad3.2}  
will be achieved.


\subsubsection{Step 1 : Initialisation - proof of Lemma \ref{D_lem:Iloc}}
\label{D_sec_Iloc}

Thanks to Assumptions $(\mathbf{P}2,3)$
(recalling that $r = \eta/4$ as stated in $(\mathbf{P}3)$),
there exists $g_\wedge >0$
such that the following lower-bound hold for any $x\in \cD_L$
and $w\in B(\aS\cdot \mathbf{e_1},\  4r)$:
\begin{equation}\label{g_wedge}
	\rho_J(x)^{-1}\cdot g(x, w) \ge g_\wedge.
\end{equation}
Let $T_1^J$, $T_2^J$ be respectively
the first and second time of jump of $X$.
Thanks to Assumptions $(\mathbf{P}1,2)$,
exploiting that $\cD_L$ is convex,
there exists $p_\wedge>0$
such that the following lower-bound holds 
for any $x\in \cD_L$ and $t\ge 0$ such that $x - v\cdot t\cdot \mathbf{e_1} \in \cD_L$:
\begin{equation}\label{p_wedge}
\textstyle \PR_x(t< T^J_1\wedge \ext)
= \exp[- \int_0^t \rho_J(x - v\cdot s\cdot  \mathbf e_1) ds]
 \ge p_\wedge.
\end{equation}
On the other hand, thanks to Assumption $(\mathbf{P}1, 3)$,
there exists also $q_\wedge>0$
such that the following lower-bound holds 
for any $x\in \cD_L$ and $t\ge a/v$ such that $x - v\cdot t\cdot \mathbf{e_1} \in \cD_L$:
\begin{equation}\label{q_wedge}
	\textstyle \PR_x(T^J_1<  t\wedge \ext)
	\ge q_\wedge.
\end{equation}

Let $t_0:= \aS/v$, $\delta:=r/v$, $x_0 \in B(x_I, r)$
and $t_\ast\in [t_0, t_0+\delta]$.
Concerning the constraint $\tIl< T_{\cD_L}$,
note that the following set is part of 	 $B(x_I, 6 r)$:
\begin{equation*}
\cA := \Lbr x_0 - v\cdot s \cdot \mathbf{e_1} \pv s \le t_0+\delta\Rbr 
	\cup \Lbr x_0 - v\cdot s\cdot \mathbf{e_1}  +w 
	\pv s \le t_0+\delta\mVg w \in B(\aS\cdot \mathbf{e_1},\  4r)\Rbr.
\end{equation*}
Since in addition $x_I\in \cD_\ell= B(0, \ell\cdot \aS)$, $6r=3\dS/2 \le 2\aS$ and $L=\ell+2$,
this set is itself a subset of $\cD_L$.
Thus, by imposing at most one such jump,
with a jump effect $w \in B(\aS\cdot \mathbf{e_1},\  4r)$,
we keep the process inside of $\cD_L$.
Let us denote by $W = \Delta X_{T^J_1}$ the size of the first jump.
We therefore restricts our analysis to the following event:
\begin{equation}\label{eq_def_cJ2}
\cJ:=	\{T_1^J < \tIl < T_2^J\wedge \ext\}
	\cap \{W \in B(\aS\cdot \mathbf{e_1},\  4r)\}.
\end{equation}
Both $\tIl < T_{\cD_L}$ and $X(\tIl) = x_0 - v\cdot \tIl\cdot \mathbf{e_1} + W$
hold a.s. on $\cJ$. It implies
the following lower-bound for any real-valued positive 
test function $f$ on $\cX$:
\begin{equation}\label{ineq_cJ}
	\E_{x_0}\lp f[X(\tIl)]\pv \tIl < T_{\cD_L}\wedge \ext\rp
	\ge \E_{x_0}\lp f[x_0 - v\cdot \tIl\cdot \mathbf{e_1} + W]\pv \cJ\rp.
\end{equation}
Let $Y := X(T^J_1-)$ be the prosition of the process just before the first jump.
Given the definition of $p_\wedge$ in \eqref{p_wedge}
and since $\cA\subset \cD_L$,
the following lower-bound holds a.s.
on the event $\{T_1^J < \tIl \wedge \ext\}
\cap \{W \in B(\aS\cdot \mathbf{e_1},\  4r)\}$:
\begin{equation}\label{eq_min_pw}
	\textstyle \PR_Y(\tIl - T^J_1< \wtd T^J_1\wedge \wtd \ext)
	\ge p_\wedge.
\end{equation}
In combination with \eqref{ineq_cJ}
and the Markov inequality,
this implies the next lower-bound:
\begin{multline}\label{ineq_cJ2}
	\E_{x_0}\lp f[X(\tIl)]\pv \tIl < T_{\cD_L}\wedge \ext\rp
\\	\ge p_\wedge 
	\E_{x_0}\lp f[x_0 - v\cdot \tIl\cdot \mathbf{e_1}+ W ]
	\pv T_1^J < \tIl \wedge \ext, W \in B(\aS\cdot \mathbf{e_1},\  4r)\rp.
\end{multline}
Conditionally on $Y$ (or equivalently on the time $T^J_1$),
the law of $W$ is given by the measure 
$\rho_J(Y)^{-1}\cdot g(Y, w) dw$.
Given the definition of $g_\wedge$ in \eqref{g_wedge},
the law of $W$
is lower-bounded by the measure $g_\wedge \cdot \idg{B(\aS\cdot \mathbf{e_1},\  4r)}(w) dw$.
Given also the definition of $q_\wedge$ in \eqref{q_wedge},
we can derive from \eqref{ineq_cJ2} the next lower-bound:
\begin{equation}\label{ineq_cJ3}
	\E_{x_0}\lp f[X(\tIl)]\pv \tIl < T_{\cD_L}\wedge \ext\rp
	\\	\ge p_\wedge\cdot g_\wedge\cdot q_\wedge
	\cdot 	\int_{B(\aS\cdot \mathbf{e_1},\  4r)} f[x_0- v\cdot \tIl\cdot \mathbf{e_1}+ w ] dw.
\end{equation}
Recalling that $x_0 \in B(x_I, r)$ and $\tIl\in [a/v, (a+r)/v]$,
we note that the following inclusions
is valid for any $x_F \in B(x_I, r)$:
\begin{equation*}
	B(x_F, r)
	\subset B(x_I, 2 r)
 \subset	x_0- v\cdot \tIl\cdot \mathbf{e_1}+ B(\aS\cdot \mathbf{e_1},\  4r).
\end{equation*}
With the fact that \eqref{ineq_cJ3} holds for any positive test-function $f$, this concludes the proof of Lemma~\ref{D_lem:Iloc}
where $c_0 := p_\wedge\cdot g_\wedge\cdot q_\wedge>0$.
\epf

\subsubsection{Step 2 : Expansion - proof of Lemma \ref{D_lem:Rloc}}
\label{sec_lem_Rloc}
By definition, the fact that $(t, c, \hat{x})\in \cR^{(L)}$
means that the following lower-bound holds 
uniformly in $x_0\in B(x_I, r)$:
\begin{equation}\label{exp1}
\PR_{x_0} \lp X_t\in dx_1
\pv t < T_{\cD_L}\wedge \ext \rp
\ge c\; \idg{B(\hat x, r)}(x_1)\,  dx_1.
\end{equation}
We see in the above proof of Lemma~\ref{D_lem:Iloc}
that the definitions of $t_0, \delta, c_0$
can be stated in terms of $L$
uniformly in $x_I$ within $\cD_\ell$.
In particular, this implies the following lower-bound
 for any $\hat{x}\in \cD_\ell$ 
and any $x_1, x_F \in B(\hat{x}, r)$:
\begin{equation}\label{exp2}
\PR_{x_1}(X(t_0)\in dx\pv t_0 < T_{\cD_L}\wedge \ext)
\ge c_0\; \idg{B(x_F, r)}(x)\, dx.
\end{equation}
In combination with the Markov property,
\eqref{exp1} and \eqref{exp2}
imply the following lower-bound 
uniformly in $x_0\in B(x_I, r)$:
\begin{equation*}
	\PR_{x_0} \lp X(t+t_0)\in dx
	\pv t+t_0 < T_{\cD_L}\wedge \ext \rp
	\ge c\cdot c_0\cdot Leb(B(\hat{x}, r))\cdot\idg{B(x_F, r)}(x)\,  dx.
\end{equation*}
Note that $c_a = c_0\cdot Leb(B(\hat{x}, r)) = c_0\cdot r^d\cdot Leb(B(0, 1))$
is a postive constant independent of $\hat{x}$.
With $t_a = t_0$,
this implies that $\{t+t_a\}\ltm \{c\cdot c_a\}\ltm B(\hat x, r)\subset \cR^{(L)}$, 
concluding the proof of Lemma \ref{D_lem:Rloc}.
\epf

\subsubsection{Step 3: Time-adjustment- proof of Lemma \ref{D_lem:PropA}}
\label{sec_expand}

Exploiting the Markov property 
as in the previous step,
the properties of $(t,c, x_I)$ and $(t', c', x')$
imply the following lower-bound 
for any $L\ge 3$, $x_I\in \cD_\ell$ and $x_0 \in B(x_I, r)$:
\begin{equation*}
	\PR_{x_0} \lp X(t+t')\in dx
	\pv t+t' < T_{\cD_L}\wedge \ext \rp
	\ge c\cdot c'\cdot Leb(B(x_I, r))\cdot\idg{B(x', r)}(x)\,  dx.
\end{equation*}
With $c_P := r^d\cdot Leb(B(0, 1))>0$, 
this concludes the proof of Lemma \ref{D_lem:PropA}.
\epf

\subsubsection{Step 4: Conclude the proof of Proposition \ref{D_lem.Ad3.2}}

Let us first consider $x_I$ and $x_F$ any points of $\cD_\ell$.
Let $K > \Lfl \|x_F - x_I\| /r\Rfl$.
For $0\le k\le K$, let \mbox{$x_k:= x_I + k\,(x_F-x_I)/K$,}
which belongs to $\cD_\ell$ (by convexity of the set).
Exploiting  Lemma \ref{D_lem:Iloc},
we choose $(t_0, c_0)$ such that
$(t_0, c_0, x_I)  \in \cR^{(L)}$.
Since $x_{k+1} \in B(x_k, r)$ for each $k$,
 by induction for any $k\le K$
thanks to Lemma \ref{D_lem:Rloc}, 
there exists $t_k, c_k >0$ such that:
$(t_k, c_k, x_k) \in \cR^{(L)}$.
In particular with $k=K$, 
there exists $t_f, c_f>0$ such that $(t_f, c_F, x_F) \in \cR^{(L)}.$

A priori, these constants $t_f, c_f>0$ still depend on the choices of $x_I$ and $x_F$.
We thus look at adjusting the different values $t_f, c_f$ 
with a finite number of such focal trajectories.
By compactness, there indeed exists $(x^j)_{j \le J}$ such that:
$\cD_\ell \subset \medcup_{j \le J} B(x^j, r)$.
Let $t_{\veebar}$ be the larger time $t_f$
needed to reach 
the vicinity of any $x_F \in \{x^j\}$ 
from any $x_I \in \{x^{j'}\}$.
 
To adjust the arrival time, 
we make the process stay some time around $x_I$.
Thanks to Lemmas \ref{D_lem:Iloc} and  \ref{D_lem:PropA},
we deduce iteratively in $k\ge 1$ the following inclusion $[k\, t_0, k\, t_0 + k\, \delta]\times \{x_I\} 
\subset \cR^{(L)}(c^k)$.
The time-intervals $[k\, t_0, k\, t_0 + k\, \delta]$ for $k\ge 1$
cover $[t_A, \infty)$ 
where $t_{A}:= t_0 \cdot \Lfl 1+ t_0/\delta\Rfl$.
Thus, for any $t \ge t_{A}$, we can find some $c>0$ for which 
$(t, c, x_I) \in \cR^{(L)}$.
Thanks to Lemma \ref{D_lem:PropA}, it
ensures, with $t\iMix:= t_{\veebar} + t_{A}$,
that there exists $c\iMix>0$
such that the following lower-bound on the density
at time $t\iMix$ holds for any $j, j' \le J$
and any $x_0 \in B(x^{j'}, r)$:
\begin{equation*}
 \PR_{x_0} \lc X(t\iMix)\in dx
 \pv t\iMix  < T_{\cD_L} \wedge \ext\rc
 \ge c\iMix\; \idg{B(x^j, r)}(x)\,  dx.
\end{equation*}
Since  $\cD_\ell \subset \medcup_{j \le J} B(x^j, r)$,
this completes the proof of Proposition \ref{D_lem.Ad3.2}.
\epf

Now that Propositions \ref{D_prop:AEF} and \ref{D_lem.Ad3.2}
are proved, as mentioned just after their statements,
the proof of Theorem  \ref{D_th.A3} is achieved.
As noted after the statement of Theorem  \ref{D_th.A3},
this also ends the proof of Theorem~\ref{D_prop.A3}.

  \section{The case of jumps occurring as in a Gibbs sampler}
  \label{D_sec_Ad5}
  
  \subsection{The core typical example}


$X$ is a pure jump process on $\bR^d$, 
for $d\ge 2$, 
again confronted to a death rate at state $x$ given by $\rho_e(x):= \Ninf{x}^2$, 
where $\Ninf{x}:= \Tsup{i\le d} |x_i|$. 
Jumps are restricted to happen 
along the vectors of an orthonormal basis 
$(\mathbf{e}_1,..., \mathbf{e}_d)$.
Independently of these directions
and of previous jumps, 
each jump occurs at rate 1
and its size follows an exponential distribution with mean $\sigma$.
This entails the following representation:
\begin{equation}
X_t 
:= x +  \Tsum{i\le N_t} \sigma W_i\, \textbf{e}_{D_i}.
\label{D_EAsy}
\end{equation}
In this formula,
 $x$ is the initial condition,
 $N_t$ a standard Poisson process on $\bZ_+$,
 while, for any $i\ge 1$,
  $W_i$ is a standard normal random variable on $\bR$,
 $D_i$ is uniform over $\II{1, d}$. 
 Moreover, all these random variables are independent from each others.

\begin{theo}
	\label{D_prop.Ad5}
		Consider $P$ the semi-group associated to the process $X$ 
	given by equation \eqref{D_EAsy}
	and weighted by the extinction event at rate $\rho_e$.
	Assume that $\sigma \le 1/8$.
Then, $P$ displays a uniform exponential quasi-stationary convergence
	with some characteristics
	$(\alpha, h, \lambda) \in \M_1(\bR^d)\times B(\bR^d)\times \bR_+$
	(cf Definition \ref{UEQS}).
	Moreover, $h$ is positive an bounded.
\end{theo}

\subsection{The main required properties}

The process under consideration in Section \ref{D_sec_Ad5}
is a specific instance of pure jump processes.
We refer to \cite{CLW17}
for a detailed presentation 
of existence results of a QSD for a pure jump process.
Contrary to the former approaches 
given in \cite{HNV94}, \cite{C10}, \cite{Co13}, \cite{CDM13}, \cite{SX15} or \cite{Sm14}
and relying on adaptations of the Krein-Rutman theorem,
the one in \cite{CLW17} exploits some maximum principle,
which makes it possible to obtain uniqueness.
It appears that no quantitative results of convergence are known.

To our knowledge, 
the restriction of having jumps only along specific directions
seems  not to have been analyzed until the current article.
As a motivation , the process $X$ could 
for instance characterize an ecosystem 
where each coordinate corresponds to a single species.

\begin{rem}
	In order to prevent concentration effects,
	several assumptions are additionally provided
	by the authors, see also \cite{BCL17}.
	Their connection to our assumption of  "almost perfect harvest"
	is a topic of interest for a future work.
\end{rem}

To highlight the generality of our approach,
we specify also in this case
 the main properties 
that we exploit.
 Let $(X_\tp)_{\tp\ge 0}$ be the pure jump process
  on $\cX:= \bR^d$ defined by:
 \begin{equation}
 X_\tp:= x + \sum_{i\le d} \int_{[0,\tp] \times \bR\times \bR_+} 
 	w \,\mathbf{e_i}\, 	\idc{ u \le g_i(X_{s^-}, w)}
 	M_i(ds,\, dw,\,du),
 	\label{D_Xdim}
 \end{equation}
 where $M_i$ are mutually independent PRaMes
 on $\bR_+\times \bR\times \bR_+$
 with intensities  $ds\,dw\, du$,
 and the $(g_i)_{i\in \II{1, d}}$
 are real-valued measurable function on $\bR^d\times \bR$.
The process is also associated
  to a state-dependent extinction rate given by $\rho_e: \bR^d\mapsto \bR_+$.
  
  In our focal example, $\rho_e(x):= \Ninf{x}^2$ and $g_i$ is defined as follows for any $i\le d$:
  $$\textstyle g_i(x, w):= \dfrac{1}{\sqrt{2\pi}\sigma}
  \exp\Big( -\frac{|w|^2}{2 \sigma^2}\Big).$$
  
Remark that the infinitesimal generator $\mathscr M$ of such generic process 
  is defined on all $C^1$ and bounded function $f$ on $\bR^d$ as follows:
  \begin{equation*}
  	\textstyle
\mathscr M f(x) := \sum_{i\le d} \int_{\bR^d} (f(x+w \mathbf{e_i}) - f(x))\cdot  g_i(x, w) dw
  	- \rho_e(x) f(x).
  \end{equation*}
  
 $(X_\tp)_{\tp\ge 0}$ 
 is a Markov Process with piecewise constant trajectories. 
 Conditionally upon $X_\tp =x$, 
 the waiting-time and size of the next jump are independent, 
 the law of the waiting-time is exponential of rate 
 $\rho_J(x):= \sum_{i\le d} \rho_J^i(x),$
  where the jump rate can be decomposed along each direction $i\in \II{1, d}$
  according to $\rho_J^i(x) 
:= \int_{\bR^d} g_i(x, w)\ dw
  <\infty$.
The jump occurs on the $i$-th coordinate 
with probability $\rho_J^i(x) / \rho_J(x)$,
then with size given by $g_i(x, w) / \rho_J^i(x)\, dw$.

 \begin{itemize}
 \item \textbf{Assumption $(\mathbf{J})$} (for Jumps)
 \item[$(\mathbf{J}1)$] The global jump rate $\rho_J$
 is upper-bounded locally in $x$.

 \item[$(\mathbf{J}2)$] Locally in $x$,
 there exists $\dS>0$
 such that the restriction of $g$ 
 to $\cX\times B(0, \dS)$
 is  lower-bounded.
 
 \item[$(\mathbf{J}3)$] The jump size has a tight law locally in $x$.
 
 \item[$(\mathbf{J}4)$] The density for each jump vector
 is  upper-bounded locally in $x$.

\item[$(\mathbf{J}5)$] 
The probability that each direction gets involved in the jump
is uniformly lower-bounded.

 \item[$(\mathbf{J}6)$]
   $\rho_e$ is bounded away from zero by $\rho > \rho\iSv$ outside some compact set.
  Moreover,  $\rho_e$ is locally bounded 
  and explosion implies extinction:
  $\ext \le \Tsup{\ell\ge 1} T_{\cD_\ell}$.

  \item[$(\mathbf{J}7)$] No stable subset: 
  $\rho\iSv < \rho\iFx$ where $\rho\iFx :=  \inf_{\{x\in \bR^d,\, i\le d\}} 
  \Lbr \rho_J^i(x)+\rho_e(x) \Rbr.$
 \end{itemize}

 \begin{theo}
 	\label{D_th.Ad5}
Provided the above conditions  $(\mathbf{J})$ are satisfied,
$P$ displays a uniform exponential quasi-stationary convergence
with some characteristics
$(\alpha, h, \lambda) \in \M_1(\bR^d)\times B(\bR^d)\times \bR_+$.
Moreover, $h$ is positive and bounded.

Besides, the Q-process exists 
and is exponentially ergodic with weight $1/h$
as stated in Corollary \ref{D_QECV},
while the uniformity in the localization procedure
holds as stated in Theorem~\ref{D_Approx} 
for $\cD_{\ell}:= B(0, \ell)$.
 \end{theo}
 We also refer to Subsections \ref{D_CRD}
for the connection with reaction-diffusion equations,
which holds in the same way, 
this time for a non-local dispersion operator of the form:
\begin{equation*}
	\mathscr M^\star u(x) 
	:=  \sum_{i\le d} 
	\lc\int_{\bR}\,  g_i(x - w_i\, \mathbf{e_i},\, w_i)
	\, u(x - w_i\, \mathbf{e_i})\, dw_i
	- \lp \int_{\bR}\,  g_i(x, w_i)\,dw_i\rp\,  u(x)\rc
\end{equation*}

Theorem \ref{D_prop.Ad5} is deduced from  Theorem \ref{D_th.Ad5}
once we prove that the process indeed satisfies $(\mathbf{J})$,
for which only $(\mathbf J7)$ is not elementary.
Let us first clarify the meaning of these assumptions.
\\

The main difference with Assumption $(\mathbf{P})$
is in the last three assumptions.
Remark that we could not avoid the comparison 
with the quantity $\rho\iSv$ in the last two. 
Note however that,
as we can see in Subsection \ref{sec_proof_J7},
Lemma 3.0.2
of \cite{AV_QSD} provide an efficient way to get upper-bounds of $\rho\iSv$.
It is hopefully enough to ensure such properties as $(\mathbf{J}6, 7)$.
While $(\mathbf{J}6)$ plays a similar role as $(\mathbf{P}1)$,
$(\mathbf{J}5, 7)$ are really specific to the fact that some jump directions are restricted.
In this illustrative example,
$(\mathbf{J}5)$ is required to make sure 
that the different directions can efficiently be explored by the process.
Thanks to it,
there exists $p_\wedge$ such that
the following lower-bound holds uniformly for $x\in \bR^d$ and  $i\le d$:
$\rho_J^i(x)  \ge  p_\wedge \cdot \rho_J(x).$ 

$(\mathbf{J}7)$ is required to prevent some directions in $\bR^d$ 
from being avoided by the process,
meaning that the probability to do so for a long time becomes negligible
even compared to extinction.

To be clear with the other properties, 
let us consider a compact $K$ subset of $\cX$.
Thanks to $(\mathbf{J}1)$,
there exists $\roV>0$ such that the following upper-bound holds uniformly for $x\in K$:
$\rho_J(x) \le  \roV$.
This is the same as $(\mathbf{P}2)$.

Thanks to $(\mathbf{J}2)$,
there exist $r, g_\wedge>0$
such that the following upper-bound holds uniformly for $x\in K$,
$i\in \II{1, d}$ and $w\in B(0, r)$: $g_i(x, w)\ge  g_\wedge$.
This is analogous to $(\mathbf{P}3)$ except that there is no drift to compensate here.

Thanks to $(\mathbf{J}3)$,
for any $\eps >0$, there exists $w_{\vee}$ such that
the following upper-bound holds uniformly for $x\in K$,
and $i\in \II{1, d}$:
$\int_{\bR} g_i(x, w) \,\idc{\Ninf{w}\ge w_{\vee}}\, dw
\le  \eps\cdot \rho_J(x)$.
This is analogous to $(\mathbf{P}4)$. 

Finally, thanks to  $(\mathbf{J}4)$
there exists  $g_\vee$ such that
the following upper-bound holds uniformly for $x\in K$, 
$i\in \II{1, d}$ and $w\in\bR$: $g_i(x,w) \le g_\vee\cdot  \rho_J(x)$.
This is the same as $(\mathbf{P}5)$.

 \subsection{Proof of $(\mathbf{J}7)$\, for our typical example}
 \label{sec_proof_J7}
 In this example, $\rho^i_J \equiv 1$ 
 while the minimal value of $\rho_e$ is simply 0. 
 Thus, we need to prove that provided $\sigma \le 1/8$, 
 $\rho\iSv < \rho_F = 1$ holds.
We rely on the criteria proposed in Lemma 3.0.2
 of \cite{AV_QSD}
 and aim at finding some set $\cD\iSv$, 
 $L\ge 1$ 
 and $t>0$ such that:
 \begin{equation*}
 \Tinf{x\in \cD\iSv} \PR_x(X_t \in \cD\iSv \mVg t < \ext \wedge T_{\cD_L})
 > e^{- t}.
 \end{equation*}
 
 We justify next our choice of $t:= (4/3)\cdot d\cdot \log 4$,
 $\cD\iSv$ and $\cD_L$
being  of the form respectively $B(0, a)$ and $\bar{B}(0, 2 a)$ for $a:= 1/4$.
 Since the jumps of $X$ occur 
 at a uniform rate $1/d$ along each direction
 and with a distribution independent of the position,
 the increase process $(X^i_t - x^i)_{t\ge 0, i\le d}$ 
 on each coordinate can be expressed
 as a i.i.d. family of processes whose law is given by:
$$
 Y_t:= \sigma \Tsum{j\le N'_t}   W_j,$$
 where $(N'_t)_{t\ge 0}$ a standard Poisson process on $\bZ_+$
 with intensity $(1/d)$
 while for any $j\ge 1$
 $W_j$ is an normal random variable.
 $N'$ and the family $(W_j)_j$ are independent.
 We remark that $Y$ is a martingale 
 with predictable quadratic variation $\LAg Y \RAg_t:= \sigma^2 t/d$,
 with the same law as $-Y$ as a symmetry.
 
Exploiting also the fact that $\rho_e$ is upper-bounded by 
$4 a^2 = 1/4$ on $\cD_L$
and thanks to the symmetries of the process,
we deduce:
 \begin{align*}
 &\Tinf{x\in \cD\iSv} \PR_x(X_t \in \cD\iSv \mVg t < \ext \wedge T_{\cD_L})
 \\&\hcm{0.5}
 \ge e^{- t / 4}\, 
 \Tinf{x\in \cD\iSv} \PR_x(\Tsup{s\le t} \Ninf{X_s - x}\le a\mVg
 \frl{i\le d}
x_i\cdot (X^i_t - x_i) \le 0)
 \\&\hcm{0.5}
  \ge e^{- t / 4}\, \lc (1/2)\cdot \PR(\Tsup{s\le t}|Y_s| \le a) \rc^d.
 \end{align*}
 Thanks to the Doob inequality, and recalling our expressions for $a$ and $t$:
 \begin{equation*}
 \PR(\Tsup{s\le t}|Y_s| \ge a)
 \le \dfrac{\E[\LAg Y \RAg_t]}{a^2} = \dfrac{16 \sigma^2 t}{d}
= \dfrac{64 \sigma^2 \log 4}{3}
< 1 / 2,
 \end{equation*}
 provided $\sigma \le  1/8 < \sqrt{3/[128\cdot \log(4)]}$.
 Since the definition of $t$ is made such that \\
 $e^{-t/4}/4^d \ge e^{-t}$,
 this concludes the following uniform lower-bound:
  \begin{equation*}
 \Tinf{x\in \cD\iSv} \PR_x(X_t \in \cD\iSv \mVg t < \ext \wedge T_{\cD_L})
 > e^{- t}.
 \end{equation*}
 Thanks to Lemma 3.0.2  of \cite{AV_QSD},
$\rho_S$ is thus necessarily smaller than $1$, 
which concludes the proof of $(J7)$ for our example.
\epf
 
  \begin{rem}
 i) 
  The condition $\sigma \le 1/8$ 
  comes only from the way we prove $(J7)$
  and is likely not to be optimal. 
  For too large values of $\sigma$ however, 
  singular concentration effects around 0 
  may play a substantial role,
  as in \cite{BCL17}.
  The event consisting 
  of forbidding any jump when starting at a Dirac Mass around 0
 might lead to a lower rate of decay in probability
   than the one consisting 
   of accumulating jumps,
  because these jumps mostly send the process to deadly regions.
  	
  ii) The purpose of assumption $(\mathbf{J}7)$\,
  is to bound the time $T_c$ at which
  either extinction occurs or all of the coordinates have changed.
  Assumption $(\mathbf{J}7)$\, indeed ensures
  an exponential moment with parameter $\Rsv$
  (cf \eqref{D_EF} below). 
  \end{rem}

 \subsection
 [Proof of the exponential convergence]
 {Proof of Theorem \ref{D_th.Ad5}}
 \label{D_sec_th5}

 For this example, 
 we consider  the family $(\cD_\ell)_{\ell\ge 1}$
as the open balls $\cD_\ell:= \bar{B}(0, \ell)$,
 now for the supremum norm $\|.\|_\infty$ for commodity.
 
 \begin{rem}
Because this norm is equivalent to the Euclidian norm, 
it is not difficult to see 
that the statements of Assumption $\mathbf{(A_F)}$
are actually equivalent for these two choices.
 \end{rem}

 Assumption \textup{$(\overline{A0})$} is clearly satisfied.
 The proof of \textup{$(A1)$} 
as stated in the following proposition
is very similar to the one of Proposition \ref{D_lem.Ad3.2}.
By these means,
we deal with each coordinate one by one 
so as to get a uniform lower-bound of the density
on a subspace of inductively increasing dimension.
The reader will be spared further details.
  \begin{prop}
  \label{D_MixPJ2}
  Assumptions $(\mathbf{J}1,2, 5)$ imply Assumption \textup{$(A1)$}, 
  with $\alc$ the uniform distribution over $\cD_1$.
  More generally, 
  for any $\ell\ge 1$,
 there exist $L>\ell$ and $\tp, \cp>0$ such that
 the follwoing inequality holds for any  $\frlq{x\in \cD_\ell}$:
    \begin{equation*}
    \PR_x \lc X(\tp)\in dy\pv
    \tp < \ext \wedge T_{\cD_L} \rc 
    	\ge \cp\, \idc{y \in \cD_{\ell}} \, dy
    \end{equation*}
  \end{prop}
Thanks to this proposition,
we know in particular that Assumption $(A1)$ holds true
for the uniform distribution over $ \cD_{1}$, i.e.:
$\zeta(dy):= \idc{y \in \cD_{1}}/Leb( \cD_{1}) \, dy$.
From Lemma 3.0.2 in \cite{AV_QSD}, 
we can (explicitly) deduce a strict upper-bound $\rho$ of $\rho_S$.
  Assumption $(A2)$ with this value of  $\rho$ 
  is clearly implied 
  by Assumption $(\mathbf{J}6)$
  for $E:= \cD_L$ where $L \ge 1$ is chosen sufficiently large.
$L$ is simply chosen
so that the extinction rate  outside of $E$ 
is larger  than $\rho$.
The proof of Assumption $(A3_F)$
for these choices
is a clear consequence of the next proposition,
whose proof
is given in the next subsection:  
 \begin{prop}
 	\label{D_lem.Ad5}
 	Assumption $(\mathbf{J})$ implies that for any $E\in \Dps$
 	and $\rho>0$, Assumption \textup{$(A3_F)$}  holds.
 \end{prop}
 With this result, 
 we can conclude that Assumption $\mathbf{(A_F)}$ holds true.
 By Theorem~\ref{D_AllPho}, Corollary \ref{D_QECV} and Theorem~\ref{D_Approx} 
 it directly implies Theorem \ref{D_th.Ad5}.
 
\subsubsection{Proof of Proposition \ref{D_lem.Ad5}  }
 We consider here three types of ``failed attempts".
 Either the process has not done all of its required jumps 
despite a very long time of observation,
 or there are too many of these jumps,
 or at least one of these jumps is too large.
 \\
 
 \noindent
  \textsl{Definition of the stopping times and time of observation}
  
  For $k\le d$, let $T_k^J$ the first time 
  at which (at least) $k$ jumps have occurred in different coordinates.
  On the event $\Lbr T_k^J < \ext\Rbr$ (for $0\le k\le d-1),$
and conditionally on $\F_{T_k^J}$,
we know from assumption $(\mathbf{J}7)$\,
 that $(T_{k+1}^J\wedge \ext) - T_k^J$
is upper-bounded by an exponential variable 
with rate parameter $\rho\iFx > \Rsv$.
Thus, with $\Rsv':= (\rho\iFx + \Rsv) / 2$,
we may define the finite quantity
$e_f := \lc 2\, \rho\iFx / (\rho\iFx - \Rsv) \rc^d$ 
and deduce the following upper-bound of the exponential moment
uniformly in $x\in \bR^d$: $\E_x \exp[ \Rsv'\cdot (T_d^J\wedge \ext) ]
\le e_f$.
Thanks to the Markov inequality, this implies the following upper-bound
on the probability that $T_d^J\wedge \ext$ takes large values:
\begin{equation}
\PR_x[T_d^J\wedge \ext > t]\cdot \exp[ \Rsv\cdot t]
\le e_f\, \exp[- (\rho\iFx - \Rsv)\cdot t/2],
\label{D_EF}
\end{equation}
which tends to 0 as $t$ tends to infinity.

 Let $\fl >0$. 
 Thanks to inequality \eqref{D_EF}, we can choose $\tZa>0$ such that:
   \begin{equation}
 \exp[\Rsv\,\tZa] \,\PR_x(\tZa< T_d^J\wedge \ext)
 \le \fl /3.
   	\label{D_tza4}
   \end{equation}
On the event $\Lbr \tZa < T_d^J\Rbr$,
we set $\Uza:= \infty$.
 This clearly implies that
 $\Uza \le \tZa$ holds a.s. on the event $\Lbr \tZa< \Tfl \Rbr$. 
 \\
 
 \noindent
\textsl{Upper-bound on the number of jumps}
 
Thanks to Assumption $(\mathbf{J}5)$, 
at each new jump, 
conditionally on the past until the previous jump,
there is a lower-bounded probability 
that a new coordinate gets altered.
The number $N_J$ of jumps before $T_d^J$ 
(on the event $\{T_d^J<\ext\}$) 
is thus upper-bounded 
by a sum of $d$ mutually independent geometric
random variables.
Therefore, we can define $n_J^\vee \ge 1$ such that
the following upper-bound holds uniformly in $x\in \bR^d$:
\begin{equation}
\PR_x(n_J^\vee \le N_J \mVg T_d^J<\ext) 
 \le \fl\, \exp[-\Rsv\,\tZa] /3.
   	\label{D_nJ4}
\end{equation}
We thus declare a failure if
 the $n_J^\vee$-th jump occurs 
 while $T_d^J$ still is not reached.
\\
 
 \noindent
\textsl{Upper-bound on the size of the jumps}
 
The crucial argument on the jump size is given by the following lemma.
\begin{lem}
\label{D_wvee}
Suppose that assumption $(\mathbf{J}3)$\  holds.
Consider any $L>0$, 
any $N\ge 1$ and any $\eps>0$.
Let $(W_i, i\ge 0)$ denote the time-ordered sequence of jump effects.
Then, there exists $w_{\vee}>0$ such that
the following lower-bound holds uniformly in the initial conditions 
$x \in \bar{B}(0, L)$:
$$
\textstyle \PR_x(\sup_{i\le N} \Ninf{W_i} \le w_{\vee}) \ge 1 - \eps.
$$
\end{lem}

\paragraph{Proof:}
The property is proved by induction over $N$, 
where one needs to adjust at each step both $\eps$ and $w_{\vee}$.
The initialization is directly implied by assumption $(\mathbf{J}3)$.
For some $N$ and $w_{\vee}^N>0$,
consider the event $\W_N(w_{\vee}^N)$ according to which 
the $N$ first jumps have a size that is upper-bounded by
$w_{\vee}^N>0$, that is $\W_N(w_{\vee}^N) := \{\sup_{i\le N} \Ninf{W_i} \le w^N_{\vee}\}$

Assume by the induction hypothesis that $w_{\vee}^N$ is chosen such that 
the following lower-bound is ensured uniformly 
for $x \in \bar{B}(0, L)$:
$\PR_x(\W_N(w_{\vee}^N)) \ge 1- \eps/2$.
On the event $\W_N(w_{\vee}^N)$, with $\Ninf{x}~\le~L$, we deduce that
 $\Ninf{X(T_N)} \le L + N\cdot w_{\vee}^N$.
Recall that $\F_{T^J_N}$ describe the information of the process up to its $N$-th jump time.
Thanks to Assumption $(\mathbf{J}3)$,
there exists $w_{\vee}^{N+1}\ge w_{\vee}^{N}$ 
such that the event $\{\Ninf{W_{N+1}}\le w_{\vee}^{N+1}\}$ 
occurs with probability greater than $1-\eps/2$
conditionally on $\F_{T^J_N}$ restricted to the event $\W_N(w_{\vee}^N)$
and uniformly on $x \in \bar{B}(0, L)$.
Note also the following inclusion:
$$\W_N(w_{\vee}^N) \cap \{\Ninf{W_{N+1}}\le w_{\vee}^{N+1}\}
\subset \W_{N+1}(w_{\vee}^{N+1}).$$
Thanks to the Markov property, the following upper-bound is then derived
for any $x \in \bar{B}(0, L)$:
\begin{align*}
\PR_x(\W_{N+1}(w_{\vee}^{N+1}))
&\ge \E_x\Big[ \PR_x(\Ninf{W_{N+1}}\le w_{\vee}^{N+1}\bv \F_{T^J_N})
\pv \W_N(w_{\vee}^N)\Big]
\\&\ge [1-\eps/2]^2 \ge 1 - \eps.
\end{align*}
The induction over $N$ then concludes the proof of the lemma.
\epf
\\

Thanks to Lemma~\ref{D_wvee},
we can choose a value $w_{\vee}>0$ such that 
$\fl\cdot \exp[-\Rsv\,\tZa]  /3$
is for any $x\in E$
an upper-bound of
the probability for the process starting from $x$
that there is a jump
before the $n_J^\vee$-th jump 
and $T_d^J\wedge \ext$ 
 with size larger than $w_{\vee}$.
We thus declare a failure if a jump larger than $w_{\vee}$ occurs.
From this we deduce that,
on the event $\{ \Tfl^1 =  \infty \}$,
 the process has stayed in $\cD_L$
for $L := \ell\iET + n_J^\vee\cdot w_{\vee}$.
 \\

On the event that at time $T_d^J < \ext$,
none of the three following conditions have been violated:

\noindent
$\quad(i)$ $T_d^J$ still has not occurred at time $\tZa$,

\noindent
or $(ii)$ the $n_J^\vee$-th jump has occurred,

\noindent
or $(iii)$ a jump of size larger than $w_\vee$ 
has occurred
 (before time $T_d^J $),\\
 we set $\Uza:= T_d^J$.
Otherwise $\Uza:= \infty$.
 \\

Given our construction 
(see \eqref{D_tza4}, \eqref{D_nJ4} 
and the above definition of $w_{\vee}$),
it is clear that:
\begin{equation*}
\Lbr\ext \wedge \tZa \le \Uza \Rbr
= \Lbr\Uza = \infty\Rbr
\quad  \text{ and }\quad
   \PR_{x} (\Uza = \infty, \,  \tZa < \ext) 
     \le \fl\, \exp(-\rho\, \tZa).
\end{equation*}

\AP{APUza}
The proof of Proposition \ref{D_lem.Ad5}  is then completed with Lemma \ref{D_MixPJ2} and
the following complementary lemma, whose proof constitutes the last step:
\begin{lem}
\label{D_UZaD}
Assume that $(\mathbf{J}3-7)$ hold,
with the preceding notations.
Then, there exists $c>0$ such that:
\begin{equation*}
 \PR_x ( X(\Uza) \in dx' 
 \pv \Uza < \ext)
 \le c\, \idc{x'\in \cD_L} dx'.
\end{equation*}
\end{lem}
 \epf

\subsubsection{Proof of Lemma \ref{D_UZaD}}

The proof is based on an induction on the coordinates 
affected by jumps in the time-interval $[0, \tZa]$.
We recall that, 
thanks to our criterion of exceptionality,
we can restrict ourselves to trajectories
where any coordinate is affected 
by at least one jump
in the time-interval $[0, \tZa]$,
while at most $n_J^\vee$ jumps have occurred in this time-interval.
We consider the sequence
of directions that the process 
follows at each successive jumps.
There is clearly
a finite number of such sequences.
In order to deduce the upper-bound 
on the density of $X(\Uza)$
presented in Lemma \ref{D_UZaD},
we merely need to prove 
the restricted versions
for any such possible sequence of directions.

Let $(i(k))$ for $k\le n_J \le n_J^\vee$ 
be a given sequence of directions in $\II{1, d}$
such that, at $k = n_J$,
 all the $d$ directions have been listed.
Let also $I(k)\in [\![1, d]\!]$ for $k\le n_J^\vee$,
be the sequence of random directions 
followed by the $n_J^\vee$ first successive jumps of $X$.
Let $U_k^J$ be the time of the $k$-th jump of $X$.

 Since in our model, 
 all directions are defined in a similar way,
 we can simplify a bit our notations without loss of generality 
 by relabeling some of the directions.
 Since we will go backwards 
 to progressively forget about the conditioning,
 we order the coordinates by the time they appear
 for the last time in the sequence $(i(k))_{k\le n_J}$.
 
It means that $i(n_J) = d$
and that,
up to the relabeling,
we exploit the unique non-decreasing function $j:\II{1, n_J}\rightarrow \II{1, d}$
such that for any $K\in \II{1, n_J}$,
$\Lbr i(k)\pv K\le k\le n_J\Rbr = [\![j(K), d]\!]$.
Let then $ K[j]$ be the largest integer $k\le n_J$
such that $j(k) \le j$.
With this definition,
it holds for any $j\in \II{1, d}$ that $i(K[j]) = j$ and that for any $k\in\II{K[j]+1, n_J}$,
 $i(k)\in \II{j+1, d}$.

\begin{rem}
 In our case, $n_J$ 
is naturally chosen as the first integer 
for which all the directions have been listed.
Yet, our induction argument
is more clearly stated 
if we do not assume this condition on $n_J$.
\end{rem}

We define the sequence $(\cA(k))_{k\le  n_J}$ of events
that encode the fact that $U^J_k$ has not reached $\ext\wedge \tZa$
and that, up to the $k$-th jump, 
the random sequence of directions coincide 
with the sequence $i$
and the size of the jumps remain uniformly bounded by $w_\vee$.
Namely, for $K\in \II{1, n_J}$:
$$\cA(K)
:= \Lbr U^J_{K} < \ext\wedge \tZa \Rbr 
\cap \Lbr \frl{k\le K-1} I(k) = i(k), \Ninf{\Delta X(U_k^J)} \le w_{\vee}\Rbr.$$
 Then, we look for a lower-bound
 that is uniform in $x\in E$
  on the following expectation that involves
 any given non-negative and measurable functions $(f_j)_{j\le d}$:
 \begin{equation*}
 	\textstyle
E^d:= \E_x\Big[ \prod_{j\le d} f_j[X^j(U^J_{n_J})]
  \pv \cA(n_J) \Big].
\end{equation*}
Define the information up to time $U^J_{n_J}$ deprived from the last jump size 
as follows:
$$\F^*_{U^J_{n_J}}
:= \sigma\big(\F_{U^J_{n_J-1}}, \{ I(n_J) = d \}\cap \{U^J_{n_J} < \ext\wedge \tZa \}\big).$$
To compute $E^d$,
we then need to compute the expectation of the following quantity:
\begin{equation}\label{D_APU1}
	\prod_{j\le d-1} f_j[X^j(U^J_{n_J-1})]
	\cdot \E_x \Big[ f_d[X^d(U^J_{n_J})] 
	\pv |\Delta X^d(U^J_{n_J})|\le w_{\vee}
	\bv \F^*_{U^J_{n_J}} \Big],
\end{equation}
restricted to the following event:
$$\cA(n_J-1)
\cap \Lbr U^J_{n_J} < \ext\wedge \tZa \Rbr 
\cap \Lbr  I(n_J) = d \Rbr.$$
Note that $X(U^J_{n_J}-) = X(U^J_{n_J-1})$ is $\F_{U^J_{n_J-1}}$-measurable,
since we consider a pure jump process. 
Thanks to the Markov property, 
the law of the next jump only depends on $x':=  X(U^J_{n_J-1})$
through the functions $(w \mapsto g_j(x', w))_{j\le d}$.
With the $\sigma$-algebra $\F^*_{U^J_{n_J}}$, 
we include the knowledge 
of the direction of the jump at time $U^J_{n_J}$,
so that only the size of this jump (possibly negative)
remains random. 
With $L := \ell\iET + n_J^\vee\cdot w_{\vee}$,
which is clearly independent of $n_J$ 
and of the particular choice of the sequence $(i(k))$,
we note the following containment property:
$$ \Ninf{ X^d(U^J_{n_J}-)}\vee \Ninf{X^d(U^J_{n_J})}
\le L.$$ 
This implies thanks to assumption $(\mathbf{J}4)$\,
that the following inequality holds a.s. on the event 
$\cA(n_J-1)
\cap \Lbr U^J_{n_J} < \ext\wedge \tZa \Rbr 
\cap \Lbr  I(n_J) = d \Rbr$:
\begin{equation}
\E_x \Big[ f_d[X^d(U^J_{n_J})] 
\pv |\Delta X^d(U^J_{n_J})|\le w_{\vee}
\bv  \F^*_{U^J_{n_J}}  \Big]
\le g_\vee  \int_{[-L, L]} f_d(x^d)\, dx^d.
\label{D_APU3}
\end{equation}
In what follows, 
the probability 
of the event $\{U^J_{n_J} < \ext\wedge \tZa \} \cap \{ I(n_J) = d \}$
is simply upper-bounded by 1.
Combining inequalities \eqref{D_APU1}, \eqref{D_APU3},
and our ordering with the definition of $K[j]$,
we deduce:
\begin{align*}
&E^d 
\le g_\vee  \int_{[-L, L]} f_d(x^d)\, dx^d
\cdot \E_x\Bigg[ \prod_{j\le d-1} f_j[X^j(U^J_{n_J-1})]
\pv \cA(n_J-1)\Big]
\\&\hcm{.5}
\le g_\vee  \int_{[-L, L]} f_d(x^d)\, dx^d
\cdot \E_x\Bigg[ \prod_{j\le d-1} f_j[X^j(U^J_{K[d-1]})]
\pv \cA(K[d-1]) \Big].
\end{align*}
Recall in  particular that $i(K[d-1]) = d-1$
and that $K[d-1] \le n_J^\vee$.
The procedure can be iterated as follows:
\begin{align*}
E^{(d-1)} 
&= \E_x\Bigg[ \prod_{j\le d-1} f_j[X^j(U^J_{K[d-1]})]
\pv \cA(K[d-1]) \Big]
\\&
\le g_\vee  \int_{[-L, L]} f_{d-1}(x^{d-1})\, dx^{d-1}
\cdot \E_x\Bigg[ \prod_{j\le d-2} f_j[X^j(U^J_{K[d-2]})]
\pv  \cA(K[d-2]) \Big],
\end{align*}
and so on until finally:
\begin{equation*}
E^d 
\le (g_\vee)^d 
\cdot \prod_{i\le d} \lp \int_{[-L, L]} f_{i}(x)\, dx \rp.
\end{equation*}

We then sum over all sequences $(i(k))$
possibly observed up to time $T_d^J$.
With the definition of the range of sequence $i$ up to step $n\ge 1$
as $\cR^i_n := \{i(k)\pv k\le n\}$,
the set of these sequences can be rigorously defined as follows:
\begin{equation*}
\Lbr (i(k))_{\{k\le n_J\}}\in \II{1, d}^{n_J}
\pv n_J \le n_J^\vee, \quad
n_J= \min\{n\ge 1; \cR^i_n = \II{1, d}\}\Rbr.
\end{equation*}
There are clearly less than $d^{n_J^\vee}$ possibilities
(there is a surjection from the set of all sequences of length $n_J^\vee$).
Since for any positive 
and measurable functions $(f_j)_{j\le d}$,
the following upper-bound is deduced uniformly for any $x\in E$:
\begin{equation*}
\E_x\lc \prod_{j=1}^d f_j[X(\Uza)] \pv \Uza < \ext\rc
\le d^{n_J^\vee}\cdot (g_\vee)^d \; 
\int_{\bar{B}(0, L)} \prod_{j=1}^d f_j(x_j) dx_1...dx_d,
\end{equation*}
it is classical that it implies
the following lower-bound on the marginal density:
\begin{equation*}
\frlq{x\in E}
\PR_x\lc X(\Uza) \in dx \pv \Uza < \ext\rc
\le d^{n_J^\vee}\cdot (g_\vee)^d \; 
\idc{x\in \bar{B}(0, L)}\; dx.
\end{equation*}
It concludes the proof of Lemma \ref{D_UZaD}.
\epf

Recall with the statement just before Lemma~\ref{D_UZaD}
that the proof of Proposition \ref{D_lem.Ad5}
is now completed.
Note also with the statement just after Proposition \ref{D_lem.Ad5}
that the proof of Theorem \ref{D_th.Ad5}
is then completed. With the statement just after Theorem \ref{D_th.Ad5},
this also concludes the proof of Theorem \ref{D_prop.Ad5}.

\section{Discussion}
\label{sec_disc}

\subsubsection{Assumption $(A3_F)$ of "Almost Perfect Harvest"}
\label{D_sec_Rafg}

\paragraph{How the different parameters have to be adjusted?}
In fact,
we will exploit this assumption only for a given single value of $\fl>0$,
which is explicitly related 
to the other parameters 
(cf Subsection~\ref{D_sec_thAf}).
But in generic proofs, 
this explicit value is not expected to be really tractable.

 The random variable $\Uza$ and $\UCa$
are thus expected to depend both on  $x\in E$ and on $\fl$, 
and to be related to $\tZa$ and $\cp$,
while these two constants must be uniform in $x$. 

\paragraph{Is it really important to consider failures?}
The purpose of introducing failure
is to handle singularities,
i.e. events 
which are rare in probability
but for which comparison estimates are poor
or simply impossible.

Notably in pure jump models,
waiting for a jump is a priori needed to loosen 
the dependency on the initial condition
(especially when the latter is a Dirac mass).
Yet, this implies that the event of a very late jump 
(being one condition of failure in the harvest)
has to be considered carefully,
to prove that its probability is negligible 
compared to the whole survival probability.
In multidimensional model, 
we may also need to wait for a jump on a specific coordinate to happen,
while there is often
a positive probability for
very singular behavior to happen meanwhile
on the other coordinates.
It is generally needed 
to adjust the singularity level
(and implicitly the efficiency of the coupling),
for the associated probability to be sufficiently small 
and for such events to be treated as a failure
in the harvest.

If this issue is made easily manageable
in the applications of Section \ref{D_sec_Adapt3}
and \ref{D_sec_Ad5},
this is mainly because we allow for both random stopping times and failure events.

Considering failures can also be of interest 
in order to exploit a Girsanov transform 
to simplify the dynamics of the process.
As can be seen in \cite{AV_Ada} and \cite{AV_M},
this transform is very efficient
to relate the original dynamics to one that is more easily described,
notably by decoupling different components of the dynamics.
	Namely,
	the original and simplified models 
	are related 
	through a change of the densities
	by a multiplicative factor
	that is upper-bounded 
	except for rare events in probability. 
	The statement of our property $(A3_F)$
	is very adapted to deal with such imprecision:
	exceptional behavior is treated as a failure, 
	so that a uniform bounds on the multiplicative factor
	can be ensured.
	From these bounds, 
	we can then deduce 
 the appropriate constant $c$ in \eqref{PtCV}.

\paragraph{Is it difficult to check the Markov properties for the harvesting time $U_{H}^\infty$?}
As in our applications, 
the definition of $U_{H}^\infty$ 
should be naturally derived from the way $\Uza$ is defined 
and the proofs that both are stopping times should be similar.

Although the law of $U_{H}^\infty$ is defined uniquely 
(which is what we need),
it is a priori unclear how to define it generally.
%
If $\Uza$ is defined 
directly in terms of the trajectories of $X$,
where $X$ has independent increments
like Brownian Motions, 
Poisson Random Measures
or say Levy processes,
$U_{H}^\infty$ can be expressed through these increments 
in the time-intervals 
$[\tau\iET^i, \tau\iET^{i+1}]$, for $i\in [\![0, \infty[\![$,
where recursively:
\begin{equation}
	\tau\iET^{i+1}:= \inf\{s\ge \tau\iET^i +\tZa; X_s \in E\}
	\wedge \ext,
	\AND
	\tau\iET^{0} = 0.
	\label{D_tdc}
\end{equation}
The Markov property on the incremental process
shall then imply the condition on 
$U_{H}^\infty-\tau\iET^1$.

\paragraph{Could we improve the assumption with less restrictions 
	on the parameters? Would it be worth it?}
The first condition in \eqref{D_FL} means indeed
that $U_H$ is required to be less than $t_F$
 for a first success in the harvest to be achieved.
This implies the following equalities:
$ \Lbr \Uza < \ext \Rbr 
=  \Lbr \Uza < \tZa \Rbr
= \Lbr \Uza = \infty \Rbr^c$.
The requirement that $\Uza$ must be less than $t_F$
is however not as stringent as it might seem
and it makes the statement of $(A3_F)$ much more tractable.

We believe it is generally compatible 
with the upper-bound on the failure event 
that one restricts any candidate $\wtd \Uza$ to be less than $t_F$,
provided $t_F$ is large enough,
possibly by reducing the considered value of $\rho$ 
 towards a value closer to $\rho\iSv$
 and by considering the event $\{\wtd\Uza\ge t_F\}$
 as an additional criterion of failure.

Nonetheless, 
a refinement of the assumption with a looser upper-bound on  $\Uza$
may still provide a better estimate of the constants involved.
Note simply that it requires 
to specify the times at which failures are stated, 
since there is no more reason for each step 
to end before time $\tZa$.
Since the statement would be much more technical,
it is not included in the current article.
Still, one may find a version of the proofs adapted for this context
in the second ArXiv version \cite{AVdisc:arXiv} of the current article.


\subsection{Brief overview of the intended applications}
\label{sec_gen}

Although greatly simplified, 
the two applications of the current article relate to eco-evolutio\-nary models.
The growth rate or the persistence of a population
is related to the individual characteristics of its members,
in other words their ``features" or ``traits". 
This effect shall be represented by the state-dependent extinction rate.
The dynamics of these traits may depend on mutations,
a changing environment or the ageing of the individuals,
for which our applications provide archetypal models.
These effects are expected to be represented 
in a ``discontinuous" fashion, i.e. with brutal transitions,
for which our approach is adapted.

Eco-evolutio\-nary models thus form a large class of applications.
Our assumption of a constant drift term
in our first application
is merely taken for simplicity
given our multidimensional state space.
Our approach could simply be adjusted
for a drift term depending on the position 
as in \cite{CH18} or \cite{CG20},
as long as it brings the process to infinity.
The proof is then much more specific 
to the biological motivation.

Such a drift term could as well be interpreted 
as the ageing of individuals
in age-structured population models
(as in \cite{TS20}),
or as the growth rate of the units
in  growth-fragmentation models (as in \cite{MMP05, BW20, CST21}).
Our hope is to see our technique fruitful for these applications
extended to a multidimensional setting
(where age or size is coupled to other individual features).
The applications are not restricted to ecology,
and may for instance come from
 chemistry (notably for polymer growth as in \cite{HY21}),
neuroscience (see the elapsed-time models e.g. in \cite{TS20})
or epidemiology (notably when the infection rate depends on the elapsed time after the infection
as presented in Subsection 1.1.2 of \cite{BP20}).

More detailed ecological models
have also been studied thanks to 
the theorems of the current paper.
Notably in \cite{AV_Ada},
we couple
a diffusive process 
specifying the population size
to a piecewise deterministic process 
specifying the adaptation of the population.
The proof is more involved than in the current paper,
notably because we use the Girsanov transform 
to decouple the diffusive 
and the piecewise deterministic components of the system.
In \cite{AV_M}, we study another related application,
in which accumulation of  deleterious mutations
is slowed down by natural selection:
the conditions of the present article are exploited
to obtain the convergence to a unique QSD
of a diffusion in an infinite dimensional state space.

More generally, 
our techniques provide 
conditions ensuring the existence and uniqueness 
of the positive eigenvector 
of general linear non-local reaction diffusion equations
(see notably Subsection \ref{D_CRD} for some partial results
and  Subsection \ref{D_MComp} for the related conditions).
The long-time behavior of structured branching processes
can typically be captured by such results,
thanks to the many-to-one formula
(see \cite{HH09}).

\subsection{Practical implications of these results}

\subsubsection{Biological motivations}

The processes presented in our applications 
can be seen as models 
for the adaptation of a population to its environment.  
In the first application, 
forcing by a regularly changing environment is considered, 
whereas in the second application, 
dependent but distinct subpopulations 
contribute to global adaptation to an otherwise fixed environment.

The environmental change in the first application
is represented by a translation 
of the fitness landscape
at constant speed $v$.
We can consider  $X$ as a summary of the individual characters of the population.
Then, the jumps come from the fixation 
of new mutations in the population,
whose rate depends on the adaptation 
of the mutant subpopulation 
(with trait $X_{t-} + w$)
as compared to the resident individuals 
(with trait $X_{t-}$).
A much more detailed description 
is proposed in  \cite{AV_Ada}.
There, 
we extend the proof to a coupled process 
involving additionally continuous fluctuations 
of the population size.

Considering distinct directions of jumps in the second application
is motivated by the interpretation that each of these directions
corresponds to the variation of a single species,
where the various $d$ species contribute to the survival 
of the community. 
Many communities are then subjects to death and reproduction events
and we can describe the state of the meta-community in this formalism
as in \cite{AV_GS}.

\subsubsection{Connection with reaction-diffusion equations}
\label{D_CRD}

The quasi-stationary regime
of the process generally prescribed 
in \eqref{D_eqS}
is expected be related to the behavior 
of the solution $(u(t, x))_{t\ge 0, x\in \bR^d}$ 
at low densities ($u \approx 0$)
to reaction-diffusion evolution equations of the form:
\begin{multline}
	\partial_t u(t, x) 
	:= v\, \partial_{x_1} u(t, x) 
	+ \int_{\bR^d}\,  g(y, x-y)\, u(t, y)\, dy
	- \lp \int_{\bR^d}\,  g(x, w)\,dw\rp\,  u(t, x)
	\\
	+ r(x, u(t, x))\, u(t, x).
	\label{D_RDeq}
\end{multline}
At low density,
the approximation of 
the growth rate $r(x, u)$
by $r_0(x):= r(x, 0)$
is usually valid, 
so as to linearize \eqref{D_RDeq}.
Looking at the linear problem
provides a criterion 
for the possibility and rate of invasion,
cf for instance \cite{BHR05} and \cite{CDM08}.

Also, if we consider $r_0(x)$
as an upper-bound 
of $r(x, u)$ for any density $u$,
the solution $\bar{u}$
to the linear problem with $r_0$
shall provide an upper-bound of $u$ 
by maximum principle approaches.
If the eigenvalue $\lambda^\star$ 
of the linear problem is negative, 
the solution $u$ is expected to asymptotically decline 
at least quicker than at rate $-\lambda^\star$.
Thus, results such as ours have
implications as a criteria for non-persistence,
as in \cite{BCV16}, \cite{BDNZ09}, \cite{BHR05}, \cite{BHR08}, \cite{BHR09}, \cite{CC98}.

In view of these interpretations,
several authors are looking
at characterizing such eigenvalue problems
when there is possibly no regular eigenvector
(see e.g. \cite{C10}, \cite{CLW17}, \cite{GR09}, 
\cite{IRS12}, \cite{CH18}).
For now, we simply conjecture that,
provided Theorem \ref{D_th.A3} applies
(with a translation of the growth rate 
by a constant to deduce an extinction rate),
they all coincide to the value prescribed with $\laZ$.

\subsubsection{Ecological relevance of the results}

In practice, the dynamics 
for which results of quasi-stationarity can be derived
usually come as an approximation. 
The relevance of the approximation is then of course at stake,
yet the quasi-stationary regime may provide insight on the conditions
of relevance.

Considering for instance our first application in Section \ref{D_sec_Mrp},
the population is certainly doomed to extinction for too strong environmental drift.
When population size strongly declines,
the estimation of individual features  
through the marginal of $X$
is not relevant.
On the other hand, 
our result of convergence for the process $X$ is not directly affected:
it holds for any value of $v$. 
As $v$ tends to $\infty$,
we shall simply have 
the asymptotic extinction rate $\lambda$ 
going to $-\infty$.
Considering the asymptotic extinction rate, notably in comparison to the convergence rate,
can nonetheless inform about the validity of the marginal of $X$ 
to capture the indivual features.

\medskip

\appendix
\section*{Appendix: Elementary facts in the absorbed setting}
\label{sec_appendix}

\paragraph{Proof of Fact \ref{fQSD}:}
Let us demonstrate that the property given in Definition \ref{UEQS} implies 
that $\alpha$ is a QSD with extinction rate $\lambda$.
For $u>0$, let us define 
$\mu_u(ds):= e^{\lambda u} \alpha P_u(ds) - \alpha(ds)$.
Then, for any $t>0$:
\begin{align*}
	\NTV{e^{\lambda t} \alpha P_t- \alpha}
	&\le \NTV{e^{\lambda t} \alpha P_t
		- e^{\lambda (t+u)} \alpha P_{t+u}} 
	+ \NTV{e^{\lambda (t+u)} \alpha P_{t+u}- \alpha}
	\\&
	\quad = \NTV{e^{\lambda t} \mu_u . P_t} 
	+ \NTV{\mu_{u+t}} 
	\\&
	\le C\, e^{-\gamma u} (e^{\lambda t}  |\!|\!|  P_t |\!|\!| + e^{-\gamma t}).
\end{align*} 
Letting $u$ tends to $\infty$ concludes the equality 
$\alpha P_t(ds) = e^{-\lambda t} \alpha(ds)$.\epf

\paragraph{Proof of Fact \ref{D_RemOpt}}

For any $u\ge 0$, $\NTV{\mu - \alpha}	\le \NTV{\mu - u\, \alpha} + |1-u|$
because \mbox{$\mu - u\, \alpha
	= \mu - u\, \alpha + (1-u)\, \alpha$}.
On the other hand \mbox{$\NTV{\mu - u\, \alpha}
	\ge |\mu(\cX) - u\, \alpha(\cX) | = |1-u|$}.
By combining these two estimates, we conclude that 
$\|\mu - u\,\alpha\|_{TV} 
\ge \|\mu - \alpha\|_{TV} / 2$.

Let $\mu$ be such that $\LAg \mu\bv 1/\heig \RAg<\infty$ 
and define the biased probability distribution:
$$\nu(dx):= \frac{1/h(x)}{\LAg \mu\bv 1/\heig \RAg} \mu(dx).$$
\noindent
Exploiting the previous inequality, we deduce 
that for any $\proj>0$:
\begin{align*}
	\|\mu - \proj\, \beta\|_{1/\heig}
	&= \LAg \mu\bv 1/\heig \RAg\cdot \NTV{\nu - (\proj/\LAg \mu\bv 1/\heig \RAg)\,  \alpha}
	\\&
	\ge \frac{\LAg \mu\bv 1/\heig \RAg}{2}\cdot \NTV{\nu - \alpha} 
	= \frac{\|\mu - \LAg \mu\bv 1/\heig \RAg\, \beta\|_{1/\heig} }{2}.
	\SQ
\end{align*} 

\noindent
\subsection*{Aknowledgment}

I am very grateful to Etienne Pardoux, my PhD supervisor, for his great support all along the redaction of this article. 
The comments of my reviewers have in addition
been particularly helpful 
to review the structure of the article 
and focus more on its specific contribution to the field. 
Pierre Picco and Pierre-André Zitt also were of great help on this aspect.
I express my sincere thanks to them. 
I would like finally to thank 
the very inspiring meetings and discussions 
brought about by the Chair “Modélisation Mathématique et Biodiversité” of VEOLIA-Ecole Polytechnique-MnHn-FX.

\end{document}